\documentclass[aap,seceqn,citesort,MSNbibl,dvips]{arximspdf}

\doi{10.1214/10-AAP684}
\volume{20}
\issue{6}
\pubyear{2010}
\firstpage{2261}
\lastpage{2294}

\makeatletter

\newtheorem{theorem}{Theorem}[section]
\newtheorem{corollary}{Corollary}[section]
\newproclaim{definition}{Definition}[section]
\newtheorem{lemma}{Lemma}[section]
\newtheorem{proposition}{Proposition}[section]
\newproclaim{Remark}{Remark}[section]
\newproclaim{Example}{Example}[section]

\makeatother

\begin{document}
\begin{frontmatter}

\title{Extreme(ly) mean(ingful): Sequential formation of a~quality group}
\runtitle{Extreme(ly) mean(ingful)}

\begin{aug}
\author[A]{\fnms{Abba M.} \snm{Krieger}\corref{}\thanksref{t1}\ead[label=e1]{krieger@wharton.upenn.edu}},
\author[B]{\fnms{Moshe} \snm{Pollak}\thanksref{t1}\ead[label=e2]{msmp@mscc.huji.ac.il}} and
\author[B]{\fnms{Ester} \snm{Samuel-Cahn}\thanksref{t2}\ead[label=e3]{scahn@mscc.huji.ac.il}
}
\runauthor{A. M. Krieger, M. Pollak and E. Samuel-Cahn}
\affiliation{University of Pennsylvania, Hebrew University and Hebrew
University}
\address[A]{A. M. Krieger\\
Department of Statistics\\
The Wharton School\\
University of Pennsylvania\\
Philadelphia, Pennsylvania 19104\\
USA\\
\printead{e1}}
\address[B]{M. Pollak\\
E. Samuel-Cahn\\
Department of Statistics\\
Hebrew University\\
Mount Scopus Campus\\
Jerusalem 91905\\
Israel\\
\printead{e2}\\
\phantom{E-mail: }\printead*{e3}}
\end{aug}

\thankstext{t1}{Supported by funds from the Marcy Bogen Chair of
Statistics at the Hebrew University of Jerusalem.}
\thankstext{t2}{Supported by the Israel Science
Foundation Grant 467/04.}

\received{\smonth{5} \syear{2009}}
\revised{\smonth{3} \syear{2010}}

%
\begin{abstract}
The present paper studies the limiting behavior of the average score
of a sequentially selected group of items or individuals, the
underlying distribution of which, $F$, belongs to the Gumbel domain
of attraction of extreme value distributions. This class contains
the Normal, Lognormal, Gamma, Weibull and many other distributions.
The selection rules are the ``better than average'' ($\beta=1$) and
the ``$\beta$-better than average'' rule, defined as follows. After the
first item is selected, another item is admitted into the group if
and only if its score is greater than $\beta$ times the average
score of those already selected. Denote by $\overline{Y}_k$ the
average of the $k$ first selected items, and by $T_k$ the time it
takes to amass them. Some of the key results obtained are: under
mild conditions, for the better than average rule, $\overline{Y}_k$
less a suitable chosen function of $\log k$ converges almost surely
to a finite random variable. When $1-F(x)=e^{-[x^\alpha+h(x)]}$,
$\alpha>0$ and $h(x)/x^\alpha\stackrel{x \rightarrow
\infty}{\longrightarrow} 0$, then $T_k$ is of approximate order
$k^2$. When $\beta>1$, the asymptotic results for $\overline{Y}_k$
are of a completely different order of magnitude. Interestingly, for
a class of distributions, $T_k$, suitably normalized, asymptotically
approaches 1, almost surely for relatively small $\beta\ge1$, in
probability for moderate sized $\beta$ and in distribution when
$\beta$ is large.
\end{abstract}

%

\end{frontmatter}

\section{Introduction and summary}\label{sec1}
Individuals are observed sequentially. The problem of whether to
accept an individual at the time that she is observed has a rich
literature. The most celebrated version is the ``Secretary Problem,''
where the criterion is to select one individual and the objective is
to maximize the probability that the best individual is chosen. This
setting has been extended in various ways including selecting a
limited number of individuals and basing the reward on the rank or
score of individual(s) selected.

Another extension that has received recent attention is to select a
group of ``quality members.'' This might occur when a team of highly
qualified professionals is assembled, for example, in an academic
department or a consulting group in a specialized area. The goal is
to find good rules for either accepting or rejecting each additional
individual into the group at the time that the individual is
observed.

One such rule that has been studied is to add a new member to the
group only if this will not decrease the average quality of the
group, termed in the literature as the ``better than average
selection rule.'' This tacitly assumes that ``quality'' is measurable.
A generalization of the rule would be to only admit a new member
whose score is say $5\%$ higher than the current average. We term
the extended rules as ``$\beta$-better than average rules.'' These
rules reduce to the better than average rule when $\beta=1$, first
considered by Preater \cite{p00}, but allows say for $\beta=1.05$ to
produce a group that is even more progressively selective than when
$\beta=1$.

The assumption that is commonly made is that the quality of the
individuals are mutually independent from a common distribution. As
the horizon, $n$, tends to infinity we study the asymptotic behavior
of the average quality of the group and the rate at which the group
grows for the $\beta$-better than average rules.


The $\beta$-better than average rules are considered in Krieger,
Pollak and Samuel-Cahn \cite{kps2}, and the present paper, which can
be read independently, can be considered its natural continuation.
Sequential selection of a ``good'' group, based only upon the relative
ranks of the observations is considered in Krieger, Pollak and
Samuel-Cahn \cite{kps}. It should be noticed that the rules
considered here can be implemented without knowledge of the
underlying distribution, though their asymptotic behavior depends
strongly on that distribution. For convenience, we assume that the
first item is always selected. However, all asymptotic results
remain correct if the selection process is adopted only after a core
group of members already exists. Also, the random variable is
assumed to be nonnegative (or the process begins with the first
nonnegative observation), because negative averages multiplied by
$\beta>1$ provide a lower level for inclusion.

Two quantities are of interest: the average quality,
$\overline{Y}_k$, of the group, after $k$ items have been retained,
and $T_k$, the time (in terms of the number of observed items) it
takes to amass a group of size $k$. Our interest is in the
asymptotics of these quantities, as $k \rightarrow\infty$. This
paper, unlike \cite{kps2}, considers $F$ belonging to the extreme
value domain of attraction of the Gumbel distribution
$\exp\{-e^{-x}\}$ only. Write $1-F(x)=\exp\{-H(x)\}$. Emphasis is
given to a subset of these distributions, which are also ``stretch
exponential'' distributions, where $H(x)=x^\alpha+h(x)$, with
$h(x)/x^\alpha\stackrel{x \rightarrow\infty}{\longrightarrow} 0$,
for all $x>x_0$, for some $x_0$, where $\alpha>0$. This class
includes the Gamma and Normal distributions as particular cases.

The ``expected overshoot,''
$f(x)=E(X-x|X>x)$, plays an essential role. Our main findings are: for
the ``better than
average'' rule ($\beta=1$), under some mild conditions, the quantity
$\overline{Y}_k-G^{-1}(\log k)$ converges a.s. to a finite random
variable where $G(x)= \int_{x_0}^x 1/f(u)\,du$. These mild conditions
are satisfied in particular by the stretch exponential distributions
with $\alpha\ge1$. It is easy to show that the functions $G(x)$
and $H(x)$ are close to each other in that $G(x)=(1+o(1))H(x)$. The
convergence of $\overline{Y}_k-G^{-1}(\log k)$ is shown in Section
\ref{sec3}, where also the convergence of the sequence of expected values and
variances of $\{\overline{Y}_k-G^{-1}(\log k)\}$ is established. The
behavior of $\overline{Y}_k$ for $\beta>1$ is very different. In
Section \ref{sec4} we show, under mild conditions, that for $\beta>1$,
$\overline{Y}_k/k^{\beta-1}$ converges a.s. to a finite positive
random variable.

The behavior of $T_k$ is discussed in Section \ref{sec5}. It is shown that
for stretch exponential random variables with $\alpha> 0$,
$\beta=1$ and every $\varepsilon>0$ one has $T_k/k^{2-\varepsilon}
\rightarrow\infty$ a.s. as $k \rightarrow\infty$, while
$T_k/k^{2+\varepsilon} \rightarrow0$ a.s. When $\alpha=1$, $T_k/k^2$
converges to a finite positive random variable. The ``standardized''
variable
\[
T_k^*=T_k\Big/\sum_{j=1}^{k-1}[1-F(\beta\overline{Y}_j)]^{-1}
\]
for $\beta\ge1$ is considered and has a very interesting behavior
for the stretch exponential with $\alpha>0$. For different values
of $\beta$ we obtain different asymptotic behavior: we show that for
$1 \le\beta< 1+1/2\alpha$ the random variable $T_k^*$ converges to
1 a.s. For $1+1/2\alpha\le\beta<1+1/\alpha$ it converges to 1 in
probability. For $\beta>1+1/\alpha$ the random variable $T^*_k$
converges in distribution to an exponential mean one distribution,
while for $\beta=1+1/\alpha$ the convergence in distribution is to a
sum of conditionally independent exponential random variables. We
conclude with Section \ref{sec6}, which contains further comments and
remarks. Section \ref{sec2} contains some preliminaries. Proofs are relegated
to the \hyperref[app]{Appendix} in order to highlight the results in the paper.

\section{Mathematical preliminaries}\label{sec2}
The
observations are denoted by $X_1,X_2,\ldots$ and are i.i.d. random
variables from a common absolutely continuous distribution~$F$.
We assume that $1-F(x)>0$ for all $x < \infty$, unless stated
otherwise.

The behavior of rules will be characterized by considering two
quantities:
\begin{itemize}
\item
$T_k=$ The number of observations inspected
until the $k$th item is retained (including that item).
\item
$\overline{Y}_k=$ The average score of the first $k$ items that are
retained.
\end{itemize}
%

The $\beta$ better than average rule is defined as follows: for
fixed $\beta$ (which is suppressed in the notation) and $T_k$
defined above as the number of items observed until the $k$th
item is selected, let $T_1=1$ and $Y_1=X_1$. Define $T_k$ and $Y_k$
inductively by
\begin{eqnarray*}
T_{k+1}&=&\inf\{i>T_k\dvtx X_i > \beta\overline{Y}_k\} ,\qquad
k=1,2,\ldots,
\\
%
%
Y_{k+1}&=& X_{T_{k+1}} ,\qquad k=1,2,\ldots.
\end{eqnarray*}
It is clear that $\overline{Y}_k$ increases in $k$ for $\beta=1$. If
$\beta>1$ we assume nonnegative $X_i$ to avoid the situation that
if $\overline{Y}_k$ is negative then the cutoff to retain an
observation becomes less stringent.

\subsection{Theorems on almost sure convergence}\label{sec21}
In this subsection, first we present two theorems, that exist in
the literature, which will be useful in proving asymptotic results
for the quantities of interest. First, we shall need the following
result, due to Robbins and Siegmund \cite{rs}, quoted as follows:
\begin{theorem}
\label{siegmund.thm}
Let $(\Omega, \mathcal{F}, P)$ be a probability
space and $\mathcal{F}_1\subset\mathcal{F}_2\subset\cdots$ a sequence of
sub-$\sigma$-algebras of $\mathcal{F}$. For each $n=1,2,\ldots,$ let
$z_n, \beta_n, \xi_n$ and $\zeta_n$ be \textit{nonnegative} $\mathcal
{F}_n$-measurable random variables such that
\[
E(z_n|\mathcal{F}_{n-1})\leq z_{n-1} (1+\beta_{n-1})+\xi_{n-1} - \zeta_{n-1}.
\]
Then $\lim_{n\rightarrow\infty} z_n$ exists and is finite and
$\sum^\infty_{n=1}\zeta_n < \infty$ a.s. on
$\{\sum^{\infty}_{n=1}\beta_n < \infty$,
$\sum^{\infty}_{n=1}\xi_n<\infty\}$.
\end{theorem}
\begin{corollary}
\label{corollary2.1} Let $z_n, \beta_n, \xi_n$ and $\zeta_n$ be
nonnegative sequences of constants such that $\sum\beta_n$ and
$\sum\xi_n$ converge, and
\[
z_n\leq z_{n-1} (1+\beta_{n-1})+\xi_{n-1} - \zeta_{n-1}.
\]
Then $\lim_{n\rightarrow\infty} z_n$ exists and is finite and
$\sum^\infty_{n=1}\zeta_n < \infty$.
\end{corollary}
\begin{pf}
This follows trivially from Theorem \ref{siegmund.thm}.
\end{pf}

We also need the following theorem that appears in Feller
\cite{feller}, page 239.
\begin{theorem}
\label{feller.thm}
Let $Q_1,Q_2, \ldots$ be independent r.v.s with
$E(Q_n)=0$, and let $S_n=\sum_{i=1}^nQ_i$. If:

(1) $b_1 <
b_2<\cdots\rightarrow\infty$ are constants
and

(2) $\sum_{n=1}^{\infty}E(Q_n^2/b_n^2) < \infty$
\[
\mbox{then } b_n^{-1}S_n\rightarrow0 \qquad\mbox{a.s. as $n
\rightarrow\infty$}.
\]
\end{theorem}

%
\subsection{Classes of distributions}\label{sec22}

Preater \cite{p00} showed that when $F$ is exponential and
$\beta=1$, $\overline{Y}_k -\log k$ converges almost surely to a
Gumbel distribution. Krieger, Pollak and Samuel-Cahn \cite{kps2}
extended this result in several ways. The asymptotic behavior of
other quantities, such as $T_k$, were obtained, values of $\beta>1$
were considered and other $F$, such as the Pareto and Beta, were
analyzed.

An interesting question is how the rules behave for other
distributions $F$. This depends on the behavior of the overshoot,
$X-a|X>a$, and its expectation $f(a)$,
%
\begin{equation}
\label{eqn23} f(a):=E(X-a|X>a).
\end{equation}

%
Let $x_F=\sup\{x\dvtx F(x)<1\}$.
\begin{definition}[(See \cite{res}, Section 1.1)]
\label{ev.def}  A distribution function
$F$, belonging to the domain of attraction of the Gumbel extreme
value distribution $\Lambda(x)=\exp\{-e^{-x}\}$, $-\infty<x<\infty$,
is called a Von Mises function (VM) if there exists $x_0$ such that
for $x_0<x<x_F$ and $r>0$
%
\begin{equation}
\label{vm.eqn}
1-F(x)=r \exp \biggl\{-\int_{x_0}^x[1/f_A(u)]\,du \biggr\}:=e^{-H(x)},
\end{equation}
where
$f_A(u)>0$, $x_0<u<x_F$ and $f_A$ is absolutely continuous on
$(x_0,x_F)$ with derivative $f_A'(u)$ and $\lim_{u \nearrow
x_F}f_A'(u)=0$.
\end{definition}

Note that Definition \ref{ev.def} is identical to definition (1.3)
in Section 1.1 of \cite{res}, except that we refer to the auxiliary
function by $f_A$ instead of $f$ to distinguish it from the
expected overshoot function. To link the two functions, we define
$g(u)=f(u)/f_A(u)$. The representation of a VM distribution function
$F$ in (\ref{vm.eqn}) is equivalent to
%
\begin{equation}
\label{vm1.eqn}
1-F(x)=r \exp \biggl\{-\int_{x_0}^x[g(u)/f(u)]\,du \biggr\}:=e^{-H(x)}.
\end{equation}
It is shown in \cite{res} that a twice differentiable VM
distribution satisfies
%
\begin{equation}
\label{vmr} \lim_{x \nearrow
x_F}\frac{F''(x)(1-F(x))}{[F'(x)]^2}=-1.
\end{equation}
It can be shown that $\lim_{u \nearrow x_F}g(u)=1$ follows from
(\ref{vmr}).

Let $G(x)$ be defined by
%
\begin{equation}
\label{geqn}
e^{-G(x)}=r \exp\biggl\{-\int_{x_0}^x [1/f(u)]\,du\biggr\}.
\end{equation}
Thus,
%
\begin{equation}
\label{overshoot} f(x)=\frac{1}{G'(x)}.
\end{equation}
Note that
%
\begin{equation}
\label{eqnn1} \frac{d}{dx}G^{-1}(x)=f(G^{-1}(x)).
\end{equation}
%

It is clear by this definition that $H(x)=(1+o(1))G(x)$ as $x \nearrow
x_F$.

We shall consider only
such VM for which $x_F=\infty$. (But see Remark \ref{Remark31}.)

%


Some of our results that will follow hold for VM distributions, but
most of the results pertain to a rich subclass. Specifically:
\begin{definition}
\label{galpha} $F$ is a generalized stretched exponential
distribution if it is VM with $H(x)=cx^\alpha+h(x)$, $h''(x)$ exists
and $c>0$, $\alpha>0$ are constants where
%
%
\begin{equation}
\label{eqn1n} \lim_{x \rightarrow\infty}\frac{h(x)}{x^\alpha}=0
\end{equation}
and
%
\begin{equation}
\label{eqn5n} \lim_{x \rightarrow\infty}\frac{h'(x)}{
x^{\alpha-1}}=0.
\end{equation}
\end{definition}

This class of distributions is denoted by $\mathcal{G}_\alpha$. By
change of variables $y=c^{1/\alpha}x$ it suffices in the sequel to
consider only $c=1$.

The reason for extending the stretched exponential by adding $h(x)$
is to include many of the classical families of distributions such
as Normal, Gamma, Lognormal and Weibull. For example, the
right-hand tail probability of the standard normal behaves like
$\phi(x)/x$ by Mills' ratio where $\phi(x)$ is the standard normal
density. Hence the standard normal belongs to $\mathcal{G}_2$ with
$h(x)=\log(x)$.
\section{Average, when $\beta=1$}\label{sec3}

In this section we consider the behavior of $\overline{Y}_k$, the
average after $k$ items are retained, using the better than average
rule. The emphasis is on the random variables that are generated
from a VM distribution. In the first subsection, we consider the
almost sure behavior, and in the ensuing subsection, results for the
expectation and variance of $\overline{Y}_k$ are presented.

Let $Z_k=Y_k-\overline{Y}_{k-1}$, the ``overshoot'' over
$\overline{Y}_{k-1}$. The results are based on the following
relationship:
%
\begin{equation}
\label{eqn1}
\overline{Y}_k=\frac{(k-1)\overline{Y}_{k-1}+Y_k}{k}=\overline
{Y}_{k-1}+\frac{Z_k}{k}=\overline{Y}_{k-1}+\frac{Z(\overline{Y}_{k-1})}{k},
\end{equation}
where $Z(a)$ is distributed like $X-a|X>a$.

The results depend on the expected overshoot $f(a)=E[Z(a)]$. We
shall use the following lemma later, that gives the expected
overshoot and squared overshoot for $F$ in $\mathcal{G}_\alpha$ for
large values of $a$. Specifically:
\begin{lemma}
\label{lemmazs} If the underlying distribution is in
$\mathcal{G}_\alpha$, $\alpha>0$ then
%
\begin{equation}
\label{lem3m} \lim_{a
\rightarrow\infty}\frac{EZ(a)}{a^{1-\alpha}/\alpha}=1
\end{equation}
and
%
\begin{equation}
\label{lem3} \lim_{a \rightarrow
\infty}\frac{EZ^2(a)}{2a^{2(1-\alpha)}/\alpha^2}=1.
\end{equation}
\end{lemma}

The proof of the results uses l'H\^{o}pital's rule on
$E(Z(a))=\int_0^\infty(1-\break F_{Z(a)}(y))\,dy$ for the expected overshoot
and $E(Z^2(a))=2\int_0^\infty y (1-F_{Z(a)}(y))\,dy$.

This result implies that $f(a)=a^{1-\alpha}/\alpha[1+o(1)]$. In some
instances we need a more refined result on the rate, that is the
$o(1)$ term, which depends on $h(x)$. An easy case, as shown in the
proof of Corollary \ref{corr31}, is when $h(x)=0$, in which case the
rate of $o(1)$ is $1/a^\alpha$.

In the more general case, we want to include $h(x)$. The point of
adding $h(x)$ is to extend our results to known distributions such
as the Normal. The role that $h(x)$ plays, is that it is small
relative to $x^\alpha$.

The following lemma provides a handle on the overshoot.
\begin{lemma}
\label{handle} If $F \sim\mathcal{G}_\alpha$, then
%
\[
f(a)=\frac{1}{H'(a)} \biggl(1+O \biggl(\frac{1}{a^\alpha} \biggr) \biggr).
\]
%
Furthermore, if $h'(x)/(x^{\alpha-\varepsilon-1})$ goes to $0$, where
$0<\varepsilon<\alpha$, we have that
%
\begin{equation}
\label{overr.eqn}
f(a)=a^{1-\alpha} \biggl[1+o \biggl(\frac{1}{a^\varepsilon} \biggr) \biggr]\Big/\alpha.
\end{equation}
\end{lemma}

These conditions on $h(x)$ and its derivatives are hardly
restrictive as the intent is for $h(x)$ to be small. In particular,
if $h(x)$ is $x^\gamma$ for $\gamma<\alpha$, then all of the above
conditions hold.

The details of the proofs throughout this and the remaining sections
of the paper appear in the \hyperref[app]{Appendix}.

\subsection{Results on almost sure convergence of the mean}\label{sec31}

The main result in this subsection is that under mild conditions
$\overline{Y}_k-G^{-1}(\log k)$ converges almost surely to a finite
random variable (Theorem \ref{evalmost.thm}). This is an extension
of the result in \cite{p00} that $\overline{Y}_k-\log k$ converges
a.s. to a Gumbel distribution
when observations are generated from an exponential
distribution. Theorem \ref{thm3n}, which is simpler than Theorem
\ref{evalmost.thm}, considers only the $\mathcal{G}_\alpha$ class of
distributions. This theorem standardizes $\overline{Y}_k$ by
dividing it by a function of $k$. Theorem \ref{evalmost.thm},
however, provides a stronger result, which for the
$\mathcal{G}_\alpha$ class of distributions is applicable when
$\alpha>1$.

The following theorem requires a slight strengthening of condition
(\ref{eqn5n}).
\begin{theorem}
\label{thm3n} If the underlying distribution function is in
$\mathcal{G}_\alpha$, where $\alpha>0$, and
%
\begin{equation}
\label{eqn5nn} \lim_{x \rightarrow\infty}
h'(x)/x^{\alpha-\varepsilon-1} =0\qquad \mbox{for some $\varepsilon>0$}
\end{equation}
and
\[
\lim_{x \rightarrow\infty}h''(x)/x^{\alpha-2}=0,
\]
then
%
\begin{equation}\quad
\label{eqn2n} \lim_{k \rightarrow\infty}\frac{\overline{Y}_k}{(\log
k)^{1/\alpha}}=\lim_{k \rightarrow\infty}
\frac{\overline{Y}_k}{G^{-1}(\log k)}=\lim_{k \rightarrow\infty}
\frac{\overline{Y}_k}{H^{-1}(\log k)}=1 \qquad\mbox{a.s.}
\end{equation}
\end{theorem}

The proof considers $S_k=(A_k-1)^2$ where
$A_k=\frac{\overline{Y}_k}{(\log k)^{1/\alpha}}$. Theorem
\ref{siegmund.thm} is used to show that $S_k$ converges almost
surely. We do not believe that the strengthening of condition
(\ref{eqn5n}) by (\ref{eqn5nn}) is necessary for the conclusion to
hold, though we use it in the proof. We know from Theorem
\ref{evalmost.thm} that it is not needed for $\alpha>1$.

The second result of this subsection, Theorem \ref{evalmost.thm},
is the stronger statement that $\overline{Y}_k-G^{-1}(\log k)$
converges a.s. to a finite random variable as $k \rightarrow
\infty$. The conditions for this result are different from those of
Theorem \ref{thm3n}, but distributions in $\mathcal{G}_\alpha$ with
$\alpha>1$ satisfy the conditions of Theorem \ref{evalmost.thm}
without the additional conditions on $h$ made in Theorem
\ref{thm3n}.



%
%
\begin{theorem}
\label{evalmost.thm} Let $\beta=1$ and $F$ be a VM distribution.
Then under conditions:

\textup{(A)} $EZ^2(a)<a^\gamma$ for
some $0<\gamma<\infty$ and all $a>a_0$ and

\textup{(B)}
$f'(a) \leq0$ for all $a \geq a_0$, for some $a_0<\infty$
\[
\overline{Y}_k-G^{-1}(\log k) \qquad\mbox{converges a.s. to a finite
random variable as $k \rightarrow\infty$}.
\]
\end{theorem}

The core of the proof is to show that $[\overline{Y}_k-G^{-1}(\log
k)]^2$ converges almost surely, using Theorem \ref{siegmund.thm}.

Conditions (A) and (B) are usually satisfied for $F$ a VM
distribution when $G(x)$ increases fast enough. In particular they
hold for $F \in\mathcal{G}_\alpha$ with $\alpha>1$. That condition
(A) holds (for all $\alpha>0$) follows from Lemma \ref{lemmazs}.
Condition (B) holds since here $f(x)=1/G'(x)=(1+o(1)) \{\alpha
x^{\alpha-1} [1+\frac{h'(x)}{\alpha x^{\alpha-1}}
] \}^{-1}$,
so from
(\ref{eqn5n}) $f(x)$ is eventually decreasing. The case $\alpha=1$
holds when $h(x)$ is increasing. If $F$ has increasing failure rate
(IFR), that is, satisfies ``new better than used,'' then condition~(B) is
satisfied.
\begin{Remark}\label{Remark31}
If $x_F < \infty$, it is easy to see that $\lim_{k \rightarrow
\infty} [\overline{Y}_k-G^{-1}(\log k)]=0$ a.s. An example $F$ of a
VM distribution with $x_F <\infty$ is $1-F(x)=e^{1/x}I(x<0)$.
\end{Remark}
%
%
\begin{corollary}
\label{corr31} Let $F \sim\mathcal{G}_\alpha$ with $\alpha\ge1$
and $h(x)=0$. Then
\[
\overline{Y}_k-\log^{1/\alpha} k \qquad\mbox{converges a.s. to a
finite random variable as $k \rightarrow\infty$}.
\]
\end{corollary}
%
%
%
%
\begin{Remark}\label{Remark32}
The conclusion of Theorem \ref{evalmost.thm} does not hold for all $F
\in\mathcal{G}_\alpha$, $\alpha>0$, thus not for all VM, for example,
$H(x)=x^{1/2}$. We omit the proof.
\end{Remark}
\subsection{Results on convergence of moments}\label{sec32}

%
\begin{theorem}
\label{thm32}
If conditions \textup{(A)} and \textup{(B)} given in Theorem \ref{evalmost.thm} hold then
there exist constants
$0<b_1, b_2, b_3<\infty$ such that
\begin{eqnarray*}
[E\overline{Y}_k-G^{-1}(\log k)] &\rightarrow& b_1,
\\
E[\overline{Y}_k-G^{-1}(\log k)]^2 &\rightarrow& b_2
\end{eqnarray*}
and hence
\[
\operatorname{Var}[\overline{Y}_k-G^{-1}(\log k)] \rightarrow b_3.
\]
\end{theorem}
\section{Average, when $\beta>1$}\label{sec4}

In this section we consider the behavior of $\overline{Y}_k$ under
the more stringent condition that an observation is retained only if
it exceeds $\beta$ times the previous average, where $\beta>1$. The
main result is that $\overline{Y}_k$ must be standardized by an
entirely different quantity, namely, $k^{\beta-1}$, in order to get
a.s. convergence. For $F \in\mathcal{G}_\alpha$ this
standardization is correct for all $\alpha>0$. The result depends on
the following relationship:
%
\begin{equation}
\label{eqn41}
\overline{Y}_k=\overline{Y}_{k-1}+\frac{(\beta-1)\overline
{Y}_{k-1}}{k}+\frac{Z(\beta\overline{Y}_{k-1})}{k}.
\end{equation}

The result concerns $B_k=\frac{\overline{Y}_k}{k^{\beta-1}}$. Let
$\mathcal{F}_k$ denote the $\sigma$-field generated by
$Y_1,\ldots,Y_k$. It follows by dividing both sides of (\ref{eqn41})
by $k^{\beta-1}$ that
%
\begin{eqnarray}
\label{eqn48n}
E(B_k|\mathcal{F}_{k-1})&=&B_{k-1} \biggl(1+O \biggl(\frac{1}{k^2} \biggr) \biggr)+E \biggl(\frac
{Z(\beta\overline{Y}_{k-1})}{k^\beta} \bigg|\mathcal{F}_{k-1}
\biggr)\nonumber\\[-8pt]\\[-8pt]
&=&B_{k-1} \biggl(1+O \biggl(\frac{1}{k^2} \biggr)
\biggr)+f(\beta\overline{Y}_{k-1})/{k^\beta}.\nonumber
\end{eqnarray}
Hence if the expected overshoot is bounded, it follows from Theorem \ref
{siegmund.thm} that
$B_k$ converges almost surely. A more refined result appears in the
next subsection followed by remarks on special cases. The section
ends with results showing that under some conditions the expected
value and variance of $B_k$ also converge.

\subsection{Almost sure convergence of the mean}\label{sec41}
We first show that $B_k$ converges almost surely under more general
conditions in the following:
\begin{theorem}
\label{thm41.thm} Assume $F$ is a VM distribution. Let
$B_k=\frac{\overline{Y}_k}{k^{\beta-1}}$ and $f(x)=E(X-x|X>x)$.

\begin{longlist}
\item If $f(x)<\frac{cx}{(\log x)^{1+\varepsilon}}$, where $c>0$ and $\varepsilon
>0$, then $B_k$ converges
a.s. to a nondegenerate positive random variable.

\item If $B_k$ converges a.s., $f$ is monotone and
$\lim_{k \rightarrow\infty}E( B_k)<\infty$ then for some constant
$x_0>0$,
%
\begin{equation}
\label{eqn45}
\int_{x_0}^\infty\frac{f(x)}{x^2}\,dx < \infty.
\end{equation}
\end{longlist}
\end{theorem}

%
\begin{Remark}\label{Remark41}
(a) Note that the sufficient condition (i) of Theorem
\ref{thm41.thm} can hold also for distributions that are not VM. An
example is the Geometric distribution.

(b)
Equation (\ref{eqn45}) does not have a $\beta$ in the expression.
Also, under the more restrictive condition of bounded expected
overshoot that was used to introduce this section (which led to an
easy proof of almost sure convergence of the desired quantity),
(\ref{eqn45}) holds.
\end{Remark}

The following is a general statement about convergence of
$\overline{Y}_k$ for the stretched exponential family of
distributions.
\begin{corollary}
\label{meanstretch} Let $F \in\mathcal{G}_\alpha$, $\alpha>0$,
$\beta>1$. Then there exists a random variable $0 < W_\beta<\infty$
such that
\[
\frac{\overline{Y}_k(\beta)}{k^{\beta-1}} \stackrel{k \rightarrow
\infty}{\longrightarrow} W_\beta\qquad\mbox{a.s.}
\]
\end{corollary}
\begin{pf}
Since $\frac{f(x)}{x^{1-\alpha}/\alpha} \rightarrow1$ it follows
that $f(x)<cx^{1-\alpha}$ for some constant \mbox{$c>0$} for all $x>x_0$
for suitable choice of $x_0$. Hence the condition in (i) of
Theorem~\ref{thm41.thm} holds.
\end{pf}
%

There exist VM distributions for which $B_k$ fails
to converge a.s. to a finite limit.
The proposition below provides a general result for when $B_k$ does
not converge to a finite limit a.s.
\begin{proposition}
\label{prop41}
Let $\Psi(a)$ be an increasing positive function of
$a$ such that
%
\begin{equation}
\label{prop41.eqn} \int_{x_0}^\infty\frac{\Psi(x)}{x^2}\,dx=\infty.
\end{equation}
%
Let $B_k=\frac{\overline{Y}_k}{k^{\beta-1}}$, $Z(a) \sim X-a|X \ge
a$ and define $Z^*(a)=Z(a)/\Psi(a)$. If
there exists a constant $a_0$ and a nonnegative random variable
$V$, not identically zero, such that for all $a \ge a_0$ $V$ is
stochastically smaller than $Z^*(a)$, then $B_k \rightarrow\infty$
a.s. as $k \rightarrow\infty$.
\end{proposition}
\begin{Example}\label{Example41}
Let $1-F_X(x)=e^{-(\log x)^2/2}$, which
is easily seen to be a VM distribution. Let $\Psi(a)=a/\log(a)$. We
shall show that the conditions (and hence the conclusions) of
Proposition \ref{prop41} hold for this example.
\end{Example}
\begin{pf*}{Proof of Proposition \ref{prop41}}
It is immediate that $\int_a^\infty\frac{\Psi(x)}{x^2}\,dx=\infty$.
Furthermore,
\begin{eqnarray*}
1-F_{Z^*(a)}(x)&=&\frac{1-F_X(a+{xa}/{\log
a})}{1-F_X(a)}=\frac{\exp\{-(\log a +\log(1+{x}/{\log
a}))^2/2\}}{\exp\{-(\log a)^2/2\}}\\
&=&\exp\biggl\{-(\log
a) \biggl(\log\biggl(1+\frac{x}{\log
a}\biggr) \biggr)-\frac{1}{2} \biggl(\log\biggl(1+\frac{x}{\log a}\biggr) \biggr)^2\biggr\}\\
&>& \exp\biggl\{-x-\frac{x^2/2}{(\log a)^2}\biggr\} > \exp\{-x-x^2/2\}
\end{eqnarray*}
for all $a > e$. Hence, if $V$ is such that
$1-F_V(x)=e^{-x-x^2/2}1(x \ge0)$ then $V$ is stochastically smaller
than $Z^*(a)$ for all $a > e$. Furthermore,
\begin{eqnarray*}
\lim_{n \rightarrow
\infty}\sum_{k=1}^n\frac{1}{k^\beta}\Psi(\gamma
k^{\beta-1})&=&\sum_{k=1}^n\frac{1}{k^\beta}\frac{\gamma
k^{\beta-1}}{\log(\gamma k^{\beta-1})}\\
&=&\gamma\sum_{k=1}^n\bigl[k\bigl(\log
\gamma+(\beta-1)\log k\bigr)\bigr]^{-1}\\
&=&\infty.
\end{eqnarray*}
\upqed\end{pf*}
%

Note that since here $f(x)=[1+o(1)]x/\log x$, it follows that one
cannot take $\varepsilon=0$ in Theorem \ref{thm41.thm}(i).


%

\subsection{Convergence of moments}\label{sec42}

We now turn to showing that the expectation of $\overline{Y}_k$
suitably normalized converges to a finite limit for all random
variables that belong to the stretch exponential.

%
We first consider $EB_k$ and $\operatorname{Var} B_k$ in a general setting.
%
\begin{theorem}
\label{lemm42} Under the following three conditions:

\textup{(a)} $\operatorname{Var} X < \infty$;

\textup{(b)} $f(a)$ is nonincreasing for $a>a_0$;

\textup{(c)} $EZ^2(a) < ca$ for some $c >0$ and $a>a_0$;

$EB_k$ and $\operatorname{Var} B_k$ converge to a finite limit.
\end{theorem}
%
%
\begin{Remark}\label{Remark42}
Condition (a) always holds for nonnegative
$X$ with $F$ a VM distribution (see Exercise 1.1.1(a) of
\cite{res}). Lemma \ref{lemmazs} implies that (c) holds for any $F
\in\mathcal{G}_\alpha$ with $\alpha\ge1/2$. Condition (b) holds
for all $F \in\mathcal{G}_\alpha$, $\alpha>1$, as well as for $X
\sim \operatorname{Exp}(1)$.
\end{Remark}

The above theorem does not apply for $F \in\mathcal{G}_\alpha$,
with $\alpha\le1$. Nevertheless, $EB_k$ converges in this case as
shown in the following theorem.
\begin{theorem}
\label{lemm41} Let $B_k=\frac{\overline{Y}_k}{k^{\beta-1}}$. If $ F
\in\mathcal{G}_\alpha$, $\alpha>0$ and $\beta>1$, then $EB_k$
converges to a finite limit.
\end{theorem}
\begin{pf}
By Theorem \ref{lemm42} we need only consider the case $\alpha\le
1$. From Le\-mma~\ref{lemmazs} it follows that for some $c$ and $k$
large enough
\begin{eqnarray*}
f(\beta\overline{Y}_{k-1})&<&c\overline{Y}{}^{1-\alpha}_{k-1}=c
[(k-1)^{\beta-1}B_{k-1} ]^{1-\alpha}=c(k-1)^{(\beta-1)(1-\alpha
)}B_{k-1}^{1-\alpha}\\
& \le&c(k-1)^{(\beta-1)(1-\alpha)}[1+B_{k-1}].
\end{eqnarray*}
Substituting this into (\ref{eqn48n}) yields
%
\begin{equation}
\label{eqnthm43n}\quad E(B_k|\mathcal{F}_{k-1}) \le
B_{k-1} \biggl[1+O \biggl(\frac{1}{k^{\min(2,1+(\beta-1)\alpha)}} \biggr) \biggr]+O \biggl(\frac
{1}{k^{1+(\beta-1)\alpha}} \biggr).
\end{equation}
Taking expectations on both sides of (\ref{eqnthm43n}) and using
Corollary \ref{corollary2.1} yields the result.
\end{pf}

\section{Time until $k$ items are kept}\label{sec5}
\subsection{Discussion of the problem}\label{sec51}

In this section we turn to the second quantity of interest, $T_k$,
the number of items that are observed until $k$ items are retained.
Unfortunately, it is generally impossible to normalize $T_k$ by a
function of $k$ and achieve almost sure convergence to a
nondegenerate random variable. Instead we consider the following
quantity:
%
\begin{equation}
\label{eqn50}
T_k^*=\frac{T_k}{\sum_{j=1}^{k-1}[1-F(\beta\overline{Y}_j)]^{-1}},
\end{equation}
which depends on the averages $\{\overline{Y}_j\}$, the expectation
of which tends to 1.

The results are
obtained for the $\mathcal{G}_\alpha$, $\alpha>0$ class of
distributions. One interesting facet of the results for $\alpha\ge
1$ is that the nature of the convergence depends on $\beta$. When
$\beta$ is relatively small, $1 \le\beta<1+\frac{1}{2\alpha}$,
then the convergence is almost sure to 1. When $\beta$ is moderate
in size, $1+\frac{1}{2\alpha} \le\beta<1+\frac{1}{\alpha}$, the
convergence is to 1, in probability. Finally, if $\beta$ is large,
$\beta\ge1+\frac{1}{\alpha}$, the convergence is in distribution
to an exponential or a sum of conditionally independent exponential
random variables with means summing up to 1.


\subsection{Almost sure convergence, when $\beta=1$}\label{sec52}
%
\begin{theorem}
\label{thm51} Let $\beta=1$ and $X_i \sim F$ where $F$ is
$\mathcal{G}_\alpha$, $\alpha>0$. Then
%
\[
T^*_k=\frac{T_k}{\sum_{j=1}^{k-1}[1-F(\overline{Y}_j)]^{-1}}\rightarrow
1 \qquad\mbox{almost surely}.
\]
\end{theorem}

The proof uses Theorem \ref{feller.thm} by conditioning on the
responses $\{Y_k\}$, letting $P_j=1-F(\overline{Y}_{j-1})$,
$b_j=\sum_{i=1}^jP_j^{-1}$ and $Q_i=T_i-T_{i-1}-P_i^{-1}$ with
$T_0=0$.

Though Theorem \ref{thm51} gives no explicit order of
magnitude of the convergence of $T_k$, in terms of $k$, we get an
idea of this magnitude in the following:
\begin{corollary}
\label{cor51} For any $\delta>0$ and $F \in\mathcal{G}_\alpha$,
$\alpha>0$, $\beta=1$
\[
\lim T_k/k^{2-\delta}=\infty\quad\mbox{and}\quad \lim T_k/k^{2+\delta}=0
\qquad\mbox{a.s.}
\]
\end{corollary}
%

For the exponential distribution $\frac{T_k}{k^2}$ converges a.s. to
a limit as
shown in \cite{kps2}.

\subsection{Asymptotic results when $\beta>1$}\label{sec53}
The focus is on $T_k^*$, the number of observations that are
observed until $k$ items are retained suitably normalized as defined
in (\ref{eqn50}).

For the sake of clarity, we consider in the continuation only $F \in
\mathcal{G}_\alpha$, $\alpha>0$, where $h(x) \equiv0$, that is,
$H(x)=x^\alpha$.
\begin{theorem}
\label{betasmall}
Let $X \sim F$ where $1-F(x)=e^{-x^\alpha}$ and $\alpha>0$. Then as
$k \rightarrow\infty$:
\begin{longlist}
\item $T_k^*
\stackrel{a.s.}{\longrightarrow} 1 \mbox{ for
$1<\beta<1+\frac{1}{2\alpha}$}$,
\item $T_k^*
\stackrel{P}{\longrightarrow} 1 \mbox{ for $1+\frac{1}{2\alpha}
\le\beta<1+\frac{1}{\alpha}$}$,
\item $ T^*_k
\stackrel{D}{\longrightarrow} \operatorname{Exp}(1)$ and
$\frac{T_k}{e^{\beta^\alpha\overline{Y}_{k-1}^\alpha}}
\stackrel{D}{\longrightarrow} \operatorname{Exp}(1)$ for
$\beta>1+\frac{1}{\alpha}$.
\end{longlist}
\end{theorem}

The result for $\beta=1+\frac{1}{\alpha}$ is of a different nature,
and hence is treated separately in Theorem \ref{thmeq}. To prove
parts (ii) and (iii) of Theorem \ref{betasmall}, we compute the
limiting generating function of $T_k$, suitably standardized, and
are able to recognize the distribution for which this limit is the
generating function. The results then follow from the Continuity
theorem. This line of reasoning is also used in proving Theorem
\ref{thmeq}.


Note that for $U \sim$ Geo$(p)$,
\[
Ee^{-tU}=\frac{1}{1+({1-e^{-t}})/({pe^{-t}})}.
\]
%
We ignore the first observation which adds one to $T_k$ (this will
have no effect on the limiting distribution). Hence the resulting
random part of $T_k$ (which we refer to as $\tilde{T}_k$ ) is the
sum of, conditionally on $\{\overline Y_j\}_{j=1}^\infty$,
independent geometric random variables where
$p_j=e^{-(\beta\overline{Y}_{j-1})^\alpha}$. We have conditionally
on $\{\overline{Y}_j\}$
\[
E \bigl(e^{-t\gamma(k)\tilde{T}_k} \bigr)=\prod_{j=2}^k \biggl[1+\frac{1-e^{-t\gamma
(k)}}{e^{-(\beta\overline{Y}_{j-1})^\alpha-t\gamma(k)}} \biggr]
^{-1},
\]
where the sequence $\gamma(k)$ is positive, and will be defined as a
function of the given~$\{\overline{Y}_k\}$, according to the need in
the proof for each particular instance, but always tends to 0. Thus
%
\begin{eqnarray}
\label{eqnm61} \log Ee^{-t\gamma(k)\tilde{T_k}}&=&-\sum_{j=2}^k
\log\biggl(1+\frac{1-e^{-t\gamma(k)}}{t\gamma(k)}t\gamma(k)e^{t\gamma
(k)+(\beta\overline{Y}_{j-1})^\alpha} \biggr)\nonumber\\[-8pt]\\[-8pt]
&=&-\sum_{j=2}^k
\log\bigl[1+\bigl(1+o_k(1)\bigr)t\gamma(k)e^{(\beta\overline{Y}_{j-1})^\alpha}
\bigr].\nonumber
\end{eqnarray}
For (ii) we let
$\gamma(k)=1/\sum_{j=1}^{k-1}e^{(\beta\overline{Y}_j)^\alpha}$ and
for (iii) we let $\gamma(k)=e^{-(\beta\overline{Y}_{k-1})^\alpha}$.

We now turn to the case where $\beta=1+1/\alpha$, so that
$\beta-1=1/\alpha$. This is the only case where conditioning on the
sequence $\{\overline{Y}_k\}$ plays a role in the limiting
distribution obtained. We know from Theorem \ref{meanstretch} that
there exists a random variable~$W$, $0<W<\infty$, such that
$\overline{Y}_k/k^{1/\alpha} \stackrel{k \rightarrow
\infty}{\longrightarrow} W$ a.s. Our result will be stated in terms
of the value of $W$.
%
\begin{theorem}
\label{thmeq} Let $X \sim F$ where $1-F(x)=e^{-x^\alpha}$, $\alpha
>0$ and $\beta=1+\frac{1}{\alpha}$. Let $W=\lim
\overline{Y}_k/k^{1/\alpha}$. Then
\[
T^*_k \stackrel{D}{\longrightarrow} \sum_{j=1}^\infty R_j
\qquad\mbox{as $k \rightarrow\infty$},
\]
where, conditionally on $W=w$, the $R_j$ are independently,
exponentially distributed with mean $\mu_j$, where
\[
\mu_j=\frac{\exp(\beta w )^\alpha-1}{\exp[j(\beta w)^\alpha]}.
\]
\end{theorem}

Note that the $\mu_j$ sum to 1.

\section{Concluding remarks}\label{sec6}
The present paper extends the results in \cite{kps2} where the
Exponential, Beta and Pareto distributions are considered in detail,
to other distributions that include the Normal, Gamma and Weibull.
The results on the special distributions considered in \cite{kps2}
are ``invertible'' in the sense that rates of convergence for
$\overline{Y}_k$ and $T_k$ imply rates of convergence for the number
of items that are kept and the average of the items kept after $n$
items are observed. The results obtained for the distributions
considered here are in general not invertible in this way.
%
%
%

Preater in \cite{p00} considered the behavior of the average of the
first $k$ items that are kept, $\overline{Y}_k$, when the
distribution generating the observations is exponential and
$\beta=1$ in the $\beta$ better than average rule. He observed that
$\overline{Y}_k-\log k$ converges a.s. and in $L_2$ to a Gumbel
distribution. The behavior of this quantity for $\beta>1$ is
markedly different. When $\beta=1$, $\overline{Y}_k -\log k$
converges a.s. When $\beta>1$, $\overline{Y}_k/k^{\beta-1}$
converges a.s. In addition, the rate when $\beta>1$ holds for many
distributions, while the amount that one subtracts from
$\overline{Y}_k$ when $\beta=1$ depends on the distribution.
%
%

There are two interesting mathematical observations. First, it is
not surprising that there should be some relationship between the
domain of attraction to which the extremal distribution of $F$
belongs and the limiting distribution of $\overline{Y}_k$, since the
$Y_k$ process will, on the average, select larger and larger items.
Preater in \cite{p00} shows that $\overline{Y}_k -\log k$ and
$\max\{X_1,\ldots,X_k\} -\log k$ have the exact same limiting Gumbel
distribution when the observations are i.i.d. from an exponential
distribution (though $\overline{Y}_k$ converges a.s. and in $L_2$
while the maximum converges only in distribution). Will the limiting
distribution of $\overline{Y}_k$, and $M_k=\max\{X_1,\ldots,X_k\}$
always agree, or at least have the same rate of convergence? From
the general theory of extreme values it follows that
%
\begin{equation}
\label{eqn61}\qquad \frac{1}{f(H^{-1}(\log k))}\bigl(M_k-H^{-1}(\log k)\bigr)
\stackrel{D}{\longrightarrow} U= \mbox{Gumbel}\qquad\mbox{as $k \rightarrow
\infty$}.
\end{equation}
This should be compared with our result for $\beta=1$ (under the
appropriate conditions of Theorem \ref{evalmost.thm}),
%
\begin{equation}
\label{eqn62}\quad \overline{Y}_k-G^{-1}(\log k) \stackrel
{\mathrm{a.s.}}{\longrightarrow} \mbox{some finite random variable}\qquad\mbox{as $k
\rightarrow\infty$}.
\end{equation}
The ``normalization'' is the same in (\ref{eqn61}) and (\ref{eqn62})
if and only if $f(x) \equiv1$, that is, if and only if the tail of the
distribution of $X$ is exponential.
%
%
%

The second interesting mathematical observation is that for the Beta
and Pareto distributions, discussed in \cite{kps2}, we get the same
kind of a.s. convergence for $T_k$, after normalization (depending
on $\beta$) for all $\beta\ge1$. In the families of distributions
considered in the present paper, different kinds of asymptotic
convergence hold for different values of $\beta$. Specifically, when
$\beta$ is relatively small, the normalized quantity converges
almost surely. When $\beta$ is in the middle range, the convergence
is in probability. For large values of $\beta$ the convergence is in
distribution.

\begin{appendix}\label{app}
\section*{Appendix}
\subsection{\texorpdfstring{Proofs for Section \protect\ref{sec3}}{Proofs for Section 3}}

\mbox{}

\begin{pf*}{Proof of Lemma \ref{lemmazs}}
Consider (\ref{lem3}). For any nonnegative random variable $Q$,
$EQ^2=2 \int_0^\infty y(1-F_Q(y))\,dy$. Thus
\[
EZ^2(a)=\frac{2\int_a^\infty(x-a)e^{-H(x)}\,dx}{e^{-H(a)}}.
\]
So
%
\setcounter{equation}{0}
\begin{eqnarray}
\label{lem4} \quad\lim_{a \rightarrow\infty}
\frac{EZ^2(a)}{2a^\delta}&=&\lim_{a \rightarrow\infty}
\frac{\int_a^\infty(x-a)e^{-(x^\alpha+h(x))}\,dx}{a^\delta
e^{-(a^\alpha+h(a))}}\nonumber\\
&\stackrel{\mathrm{l}\mbox{\fontsize{8.36}{10}\selectfont{'
H\^{o}pital}}}{=}& \lim_{a
\rightarrow\infty} \frac{-\int_a^\infty
e^{-(x^\alpha+h(x))}\,dx}{e^{-(a^\alpha+h(a))}\{\delta
a^{\delta-1}-[\alpha
a^{\alpha-1}+h'(a)]a^\delta\}}\nonumber\\[-8pt]\\[-8pt]
&=&\lim_{a
\rightarrow\infty} \frac{\int_a^\infty
e^{-(x^\alpha+h(x))}\,dx}{e^{-(a^\alpha+h(a))}\alpha
a^{\alpha+\delta-1} [1+{h'(a)}/({\alpha
a^{\alpha-1}}) ]}\nonumber\\
&=& \lim_{a \rightarrow
\infty}\frac{\int_a^\infty
e^{-(x^\alpha+h(x))}\,dx}{e^{-(a^\alpha+h(a))}\alpha
a^{\alpha+\delta-1}}\nonumber
\end{eqnarray}
by (\ref{eqn5n}).
Using l'H\^{o}pital's rule once more we get that the value in
(\ref{lem4}) equals
\begin{eqnarray*}
&&\lim_{a \rightarrow\infty} \frac{-e^{-(a^\alpha+h(a))}}{\alpha
e^{-(a^\alpha+h(a))} \{(\alpha+\delta-1)a^{\alpha+\delta-2}-[\alpha
a^{\alpha-1}+h'(a)]a^{\alpha+\delta-1} \}}\\
&&\qquad=\lim_{a \rightarrow
\infty}\frac{1}{\alpha^2a^{2(\alpha-1)+\delta}}.
\end{eqnarray*}
Thus if we take $\delta=2(1-\alpha)$ the above limit is $1/\alpha^2$
and (\ref{lem3}) follows.

%
The proof for $E(Z(a))=\int_a^\infty e^{-H(x)}\,dx/e^{-H(a)}$ follows
in a similar manner.
\end{pf*}
\begin{pf*}{Proof of Lemma \ref{handle}}
Through integration by parts
\begin{eqnarray*}
f(a)&=&\frac{\int_a^\infty
e^{-H(x)}\,dx}{e^{-H(a)}}\\
&=&\frac{\int_a^\infty
e^{-H(x)}{H'(x)}/{H'(x)}\,dx}{e^{-H(a)}}\\
&=&\frac{ [-e^{-H(x)}/H'(x) ]_a^\infty}{e^{-H(a)}}-\frac{\int_a^\infty
e^{-H(x)}{H''(x)}/{(H'(x))^2}\,dx}{e^{-H(a)}}.
\end{eqnarray*}
Note that ${H''(x)}/(H'(x))^2$ tends to 0 for a VM distribution. Now
to get the rate, consider $H(x)=x^\alpha+h(x)$ where $ \lim_{x
\rightarrow\infty} h'(x)/ x^{\alpha-1}=0$ by (\ref{eqn5n}), and
assume $\lim_{x \rightarrow\infty} h''(x)/ x^{\alpha-2}=0$. This
implies that
$\frac{H''(x)}{(H'(x))^2}=O (\frac{1}{x^\alpha} )$. Since
the first term is $1/H'(a)$, using l'H\^{o}pital's rule on the
second term yields
\[
f(a)=\frac{1}{H'(a)} \biggl(1+O \biggl(\frac{1}{a^\alpha} \biggr) \biggr).
\]
Finally, to get the rate at which $f(a)/(a^{1-\alpha}/\alpha)$ goes
to 1, we need the rate at which $h'(x)/x^{\alpha-1}$ goes to 0. If
we assume that $h'(x)/(x^{\alpha-\varepsilon-1})$ goes to $0$, where
$0<\varepsilon<\alpha$, we have (\ref{overr.eqn}).
\end{pf*}
\begin{pf*}{Proof of Theorem \ref{thm3n}}
Let $A_k=\frac{\overline{Y}_k}{(\log k)^{1/\alpha}}$ and
$S_k=(A_k-1)^2$. Note that $\frac{\log(k-1)}{\log k}=1-\frac{1}{k
\log k}+O (\frac{1}{k^2} )$, thus
\[
\biggl(\frac{\log(k-1)}{\log k} \biggr)^{1/\alpha}=1-\frac{1}{\alpha
k\log k}+O \biggl(\frac{1}{k^2} \biggr).
\]
Hence,
\begin{eqnarray*}
S_k&=& \biggl(A_{k-1} \biggl(\frac{\log(k-1)}{\log
k} \biggr)^{1/\alpha}+\frac{Z_k}{k(\log k)^{1/\alpha}}-1 \biggr)^2\\
&=& \biggl[(A_{k-1}-1)-A_{k-1} \biggl(\frac{1}{\alpha k \log
k}+O \biggl(\frac{1}{k^2} \biggr) \biggr)+\frac{Z_k}{k(\log k)^{1/\alpha}} \biggr]^2\\
&=&(A_{k-1}-1)^2+A_{k-1}^2 \biggl(\frac{1}{\alpha^2k^2(\log
k)^2}+O \biggl(\frac{1}{k^3} \biggr) \biggr)+\frac{Z_k^2}{k^2(\log
k)^{2/\alpha}}\\
&&{}+2(A_{k-1}-1) \biggl[\frac{Z_k}{k(\log
k)^{1/\alpha}}-A_{k-1} \biggl(\frac{1}{\alpha k \log
k}+O\biggl (\frac{1}{k^2} \biggr) \biggr) \biggr]\\
&&{}-2\frac{Z_k}{k(\log k)^{1/\alpha}}A_{k-1} \biggl(\frac{1}{\alpha k
\log k}+O \biggl(\frac{1}{k^2} \biggr) \biggr).
\end{eqnarray*}
%
Taking conditional expectations on both
sides, using (\ref{eqn23}), we therefore get
%
\begin{eqnarray}
\label{eqn6n}
&&
E (S_k|\mathcal{F}_{k-1} ) \nonumber\\
&&\qquad\le
S_{k-1}+2(S_{k-1}+1) \biggl(\frac{1}{\alpha k^2(\log k)^2
}+O \biggl(\frac{1}{k^3} \biggr) \biggr)\nonumber\\[-8pt]\\[-8pt]
&&\qquad\quad{}+\frac{E(Z^2(\overline{Y}_{k-1})|\mathcal{F}_{k-1})}{k^2(\log
k)^{2/\alpha}}\nonumber\\
&&\qquad\quad{}+2(A_{k-1}-1) \biggl[\frac{f(\overline{Y}_{k-1})}{k (\log
k)^{1/\alpha}}-A_{k-1} \biggl(\frac{1}{\alpha k \log
k}+O \biggl(\frac{1}{k^2} \biggr) \biggr) \biggr]\nonumber
\end{eqnarray}
since $A_{k-1}^2 \le2(S_{k-1}+1)$. The first two terms in
(\ref{eqn6n}) therefore cause no problem in the application of
Theorem \ref{siegmund.thm} to $S_k$. By Lemma \ref{lemmazs}, for all
$k$ sufficiently large
\begin{eqnarray*}
\frac{E(Z^2(\overline{Y}_{k-1})|\mathcal{F}_{k-1})}{k^2(\log
k)^{2/\alpha}}
&<&\frac{2(1+\varepsilon)}{\alpha^2}\frac{\overline{Y}_{k-1}^{2(1-\alpha)}}
{k^2(\log k)^{2/\alpha}}\\
&=&\frac{2(1+\varepsilon)}{\alpha^2}\frac{A_{k-1}^{2(1-\alpha)}(\log
(k-1)^{1/\alpha})^{2(1-\alpha)}}{k^2(\log k)^{2/\alpha}}\\
&<&\frac{2(1+\varepsilon)}{\alpha^2}\frac{A_{k-1}^{2(1-\alpha)}}{k^2(\log
k)^2}<\frac{2(1+\varepsilon)}{\alpha^2}\frac{(A_{k-1}^2+1)}{k^2(\log
k)^2}
\\
&<&\frac{2(1+\varepsilon)}{\alpha^2} \biggl[\frac{2S_{k-1}}{k^2(\log
k)^2}+\frac{3}{k^2(\log k)^2} \biggr],
\end{eqnarray*}
so the second term in the last expression is summable, and again
factoring out $S_{k-1}$ the first term is also summable.

It remains to deal with the last term in (\ref{eqn6n}).
From (\ref{overr.eqn}),
\[
f(\overline{Y}_{k-1})=\frac{\overline{Y}{}^{1-\alpha}_{k-1}
[1+o(1/\overline{Y}{}^\varepsilon_{k-1}) ]}{\alpha}.
\]
%
Thus the first term in the square brackets in (\ref{eqn6n})
satisfies
\begin{eqnarray*}
\frac{f(\overline{Y}_{k-1})}{k \log
k^{1/\alpha}}&=&\frac{\overline{Y}{}^{1-\alpha}_{k-1} [1+o(1/\overline
{Y}{}^\varepsilon_{k-1}) ]}{\alpha
k (\log k)^{1/\alpha}}\\
&=&\frac{A_{k-1}^{1-\alpha}[(\log
(k-1))^{1/\alpha}]^{1-\alpha} [1+o(1/(A_{k-1}(\log
(k-1))^{1/\alpha})^\varepsilon) ]}{\alpha k (\log
k)^{1/\alpha}}\\
&=&\frac{A_{k-1}^{1-\alpha}
[1+o (A_{k-1}^{-\varepsilon}/ (\log
k)^\omega) ]}{\alpha k \log k},
\end{eqnarray*}
where $\omega=\varepsilon/\alpha$. The last line in (\ref{eqn6n}) can
therefore be rewritten as
%
\begin{equation}
\label{eqn7n} -
\frac{2(A_{k-1}-1)A_{k-1} [1-A_{k-1}^{-\alpha}+O ({\log
k}/{k} )+A_{k-1}^{-\alpha}o ({A_{k-1}^{-\varepsilon}}/{(\log
k)^\omega} ) ]}{\alpha k \log k}.\hspace*{-32pt}
\end{equation}
We want to study when (\ref{eqn7n}) is positive for large $k$. This
depends on the term in brackets, which to simplify notation we
denote by $R(x)$, where $x=A_{k-1}$, and the dependence of $R(x)$ on
$k$ is implicit. Note that for $k$ sufficiently large $O(\log
k/k)<\delta_k \equiv1/(\log k)^\omega$.
Also note that $\nu_k\equiv A_{k-1}^{-\varepsilon}\delta_k
<\overline{Y}{}^{-\varepsilon}_{k-1} \rightarrow0$
as $k
\rightarrow\infty$ and that $\nu_k<x_0^{-\varepsilon}\delta_k$ if
$A_{k-1}>x_0>0$. Hence when $k$ is sufficiently large,
$\underline{R}(x) \le R(x) \le\overline{R}(x)$ where
%
\begin{equation}
\label{eqn8l} \underline{R}(x)=1-x^{-\alpha}-\delta_k-\nu_k
x^{-\alpha}
\end{equation}
and
%
\begin{equation}
\label{eqn8r} \overline{R}(x)=1-x^{-\alpha}+\delta_k+\nu_k
x^{-\alpha}.
\end{equation}

The aim is to show that (\ref{eqn7n}) is positive only when
$1-c\delta_k \le x \le1+c\delta_k$ for a suitably chosen constant
$0<c<\infty$. We consider two cases:

\begin{longlist}
\item Assume $x>1$. Then (\ref{eqn7n}) is positive for the values of $x$
such that \mbox{$R(x)<0$}.
Since $\underline{R}(x) \le R(x)$, the values of $x$ such that
$\underline{R}(x)<0$, or, equivalently, the values of $x$ such that
$x^\alpha\underline{R}(x)<0$ include the values of $x$ such that
$R(x)<0$. It suffices to consider
%
\begin{equation}
\label{eqn8r1}
x^\alpha-1-\delta_k x^\alpha-\nu_k<0.
\end{equation}
The set of $x$ such that (\ref{eqn8r1}) holds is equivalent to the
set of $x$ such that
\[
x< \biggl(\frac{1+\nu_k}{1-\delta_k} \biggr)^{1/\alpha} < 1+
c_1\delta_k
\]
for $k$ large for a suitably
chosen constant $0<c_1<\infty$.

\item  Assume $x <1$. Then (\ref{eqn7n}) is positive for the values of $x$
such that \mbox{$R(x)>0$}.
Since $\overline{R}(x) \ge R(x)$, the values of $x$ such that
$\overline{R}(x)>0$, or, equivalently, the values of $x$ such that
$x^\alpha\overline{R}(x)>0$ include the values of $x$ such that
$R(x)>0$. Hence, we want to consider when
%
\begin{equation}
\label{eqn8r2} x^\alpha+\delta_k x^\alpha+\nu_k >1.
\end{equation}
Since $\delta_k$ and $\nu_k$ are arbitrarily small for $k$
sufficiently large, there exists $x_0>0$ such that for
(\ref{eqn8r2}) to hold it is sufficient that $x>x_0$.
Therefore, (\ref{eqn8r2})
is equivalent to
%
\begin{equation}
\label{eqn8r4}
x> \biggl(\frac{1-\nu_k}{1+\delta_k} \biggr)^{1/\alpha}>1-c_2\delta_k
\end{equation}
for $k$ sufficiently large for a suitable chosen constant $0<c_2<\infty$.

The above analysis shows that (\ref{eqn7n}) can be bounded from
above by zero when $A_{k-1}$ is outside the interval $c \pm
\delta_k$. When it is inside, (\ref{eqn7n}) is bounded by $O(1/[k
(\log k)^\omega])$. Hence (\ref{eqn7n}) is summable.
Thus $S_k$ converges a.s. by Theorem~\ref{siegmund.thm}.
%
%
%
%
If $S_k$ converges to a value different from 0 this would lead to a
contradiction, as the sum of the terms in (\ref{eqn7n}) would go to
minus infinity, while $S_k$ is nonnegative. Hence $A_k$ tends to 1
a.s.

Note that when $H(x)=x^\alpha+h(x)$, and (\ref{eqn1n}) holds, then
necessarily $H^{-1}(x)=x^{1/\alpha}+h^*(x)$ where
$\frac{h^*(x)}{x^{1/\alpha}} \stackrel{x \rightarrow
\infty}{\longrightarrow}0$. Thus $\frac{H^{-1}(\log k)}{(\log
k)^{1/\alpha}} \stackrel{k \rightarrow\infty}{\longrightarrow} 1$
and the right-hand side of (\ref{eqn2n}) follows. Since
$G^{-1}(x)=H^{-1}(x)[1+o(x)]$, also the middle term in (\ref{eqn2n})
converges a.s. to 1.\qed
\end{longlist}
\noqed\end{pf*}
\begin{pf*}{Proof of Theorem \ref{evalmost.thm}}
The first step in the proof is to show that
$(\overline{Y}_k-G^{-1}(\log k))^2$ converges a.s. to a finite
random variable as $k \rightarrow\infty$.

Since $\overline{Y}_k \rightarrow\infty$, there will be a (possibly
random) $k_0$ such that for all $k >k_0$, everything written below
holds. Consider $k >k_0$ only. Let $c_k=G^{-1}(\log k)$. Then, by
(\ref{eqnn1}) and the boundedness of $f$,
%
\begin{eqnarray}
\label{eqn4}
c_k-c_{k-1}&=&\bigl(\log k
-\log(k-1)\bigr)[G^{-1}(u_k)]'\nonumber\\
&=&-\log\biggl(1-\frac{1}{k}\biggr)f(G^{-1}(u_k))\\
&=&\frac{f(G^{-1}(u_k))}{k}+O\biggl(\frac{1}{k^2}\biggr),\nonumber
\end{eqnarray}
where the $O(\frac{1}{k^2})$ term is positive and
%
\begin{equation}
\label{eqn5} \log(k-1) \leq u_k \leq\log k.
\end{equation}
Note that the last equality in (\ref{eqn4}) follows since $f$ is
bounded by condition (B).
Now write
\begin{eqnarray*}
(\overline{Y}_k-c_k)^2&=&\biggl[(\overline{Y}_{k-1}-c_{k-1})+\frac{Z(\overline
{Y}_{k-1})}{k}
+(c_{k-1}-c_k)\biggr]^2\\
&=&(\overline{Y}_{k-1}-c_{k-1})^2+\frac{Z^2(\overline
{Y}_{k-1})}{k^2}+(c_{k-1}-c_k)^2\\
&&{}+\frac{2Z(\overline{Y}_{k-1})}{k}(c_{k-1}-c_k)\\
&&{}+
2(\overline{Y}_{k-1}-c_{k-1}) \biggl[\frac{Z(\overline
{Y}_{k-1})}{k}+(c_{k-1}-c_k) \biggr].
\end{eqnarray*}
%
Taking conditional expectation, conditional on
$\mathcal{F}_{k-1}$, yields
%
\begin{eqnarray}\label{eqn7}\quad\qquad
E[(\overline{Y}_k-c_k)^2|\mathcal{F}_{k-1}]
&=&\underbrace{(\overline{Y}_{k-1}-c_{k-1})^2}_{\mathrm{(i)}}{}+{}\underbrace{\frac
{E[Z^2(\overline{Y}_{k-1})|\mathcal{F}_{k-1}]}{k^2}}_{\mathrm{(ii)}}\nonumber\\
&&{}+{}\underbrace{(c_{k-1}-c_k)^2}_{\mathrm{(iii)}}{}+{}\underbrace{\frac{2f(\overline
{Y}_{k-1})}{k}(c_{k-1}-c_k)}_{\mathrm{(iv)}}
\\
&&{}
+{}\underbrace{2(\overline{Y}_{k-1}-c_{k-1}) \biggl[\frac{f(\overline
{Y}_{k-1})}{k}+(c_{k-1}-c_k) \biggr]}_{\mathrm{(v)}}.\nonumber
\end{eqnarray}
We shall show that the conditions for Theorem \ref{siegmund.thm}
hold. We shall examine each term in (\ref{eqn7}) separately. We
first show that for any $\omega>0$
%
\begin{equation}
\label{eqn8}
\overline{Y}_k/k^{\omega} \rightarrow0 \qquad\mbox{as $k
\rightarrow\infty$ }\mbox{a.s}.
\end{equation}
Let $W_k(\omega)=\frac{\overline{Y}_k}{k^{\omega}}$. Then clearly
$W_k(\omega)>0$ and
\[
E[W_k(\omega)|\mathcal{F}_{k-1}]=
\biggl(\frac{k-1}{k} \biggr)^{\omega}W_{k-1}
(\omega)+\frac{f(\overline{Y}_{k-1})}{k^{\omega+1}}
<W_{k-1}(\omega)+\frac{B}{k^{\omega+1}},
\]
since $f$ is bounded (where we have denoted its bound by $B$). It
follows that $W_k(\omega)$ converges a.s. to a finite limit,
$L(\omega) \ge0$. Then also $W_k(\omega/2) \rightarrow L(\omega
/2)$ a.s. But $W_k(\omega)=W_k(\omega/2)/k^{\omega/2}$, thus the
limit must be $0$ for all $\omega$.

Now consider term (ii) of (\ref{eqn7}). By condition (A) and
(\ref{eqn8}), for all $k$ sufficiently large
%
\begin{equation}
\label{eqn9}
\frac{E[Z^2(\overline{Y}_{k-1})|\mathcal
{F}_{k-1}]}{k^2}<\frac{\overline{Y}_{k-1}^\gamma}{k^2} <
\frac{\varepsilon k^{\omega\gamma}}{k^2} \qquad\mbox{a.s.}
\end{equation}
Choose $\omega<\frac{1}{\gamma}$ and write $1-\omega
\gamma=\delta$. The rightmost expression in (\ref{eqn9}) is then
$\varepsilon/ k^{1+\delta}$, which clearly is summable.

Term (iii) is summable by (\ref{eqn4}) and the boundedness of $f$.

Term (iv) is negative, and hence causes no problem.

Term (v): note first that by (\ref{eqn4})
%
\begin{eqnarray}
\label{eqn10}\quad
\frac{f(\overline{Y}_{k-1})}{k}+(c_{k-1}-c_k)&=&\frac{f(\overline{Y}_{k-1})-
f(G^{-1}(u_k))}{k}+O\biggl(\frac{1}{k^2}\biggr)\nonumber\\[-8pt]\\[-8pt]
&=&\frac{(\overline{Y}_{k-1}-G^{-1}(u_k))}{k}f'(d_k)+O\biggl(\frac{1}{k^2}\biggr),\nonumber
\end{eqnarray}
where $d_k$ is a value between $\overline{Y}_{k-1}$ and
$G^{-1}(u_k)$. Since $G^{-1}$ is increasing, it follows from
(\ref{eqn5}) that
%
\begin{equation}
\label{eqn11} c_{k-1} \le G^{-1}(u_k) \le c_k.
\end{equation}
Consider two cases:

(a) $\overline{Y}_{k-1}-c_{k-1} \le0$. Then by
(\ref{eqn11}) also $\overline{Y}_{k-1}-G^{-1}(u_k) \le0$ and by condition~(B)
(v) is negative since the $O(1/k^2)$ term is positive.

(b) $\overline{Y}_{k-1}-c_{k-1}>0$. If also
$\overline{Y}_{k-1}-G^{-1}(u_k) \ge0$,
the previous argument goes through, except that we still must show
that
$(\overline{Y}_{k-1}-c_{k-1})I(\overline{Y}_{k-1}-c_{k-1}>0)/k^2$ is
summable. Now write
$(\overline{Y}_{k-1}-c_{k-1})I(\overline{Y}_{k-1}-c_{k-1}>0)<(\overline
{Y}_{k-1}-c_{k-1})^2+1$.
Thus
%
\begin{equation}
\label{eqn12}\quad
(\overline{Y}_{k-1}-c_{k-1})I(\overline{Y}_{k-1}-c_{k-1}>0)/k^2 \le
\frac{(\overline{Y}_{k-1}-c_{k-1})^2}{k^2}+\frac{1}{k^2}.
\end{equation}
The first term on the right-hand side of (\ref{eqn12}) can be
combined with (i) in (\ref{eqn7}), and the second is clearly
summable.

Now suppose $\overline{Y}_{k-1}-G^{-1}(u_k) < 0 <
\overline{Y}_{k-1}-c_{k-1}$. Then
%
\begin{equation}
\label{eqn13} c_{k-1} < \overline{Y}_{k-1} <G^{-1}(u_k) <c_k.
\end{equation}
Since both $|\overline{Y}_k-c_{k-1}|$ and
$|\overline{Y}_k-G^{-1}(u_k)|$ are less than $c_k-c_{k-1}$, it
follows from (\ref{eqn4}) that (v) is summable.

It follows that in all cases we can write
$E[(\overline{Y}_k-c_k)^2|\mathcal{F}_{k-1}] \le
(\overline{Y}_{k-1}-c_{k-1})^2(1+B_{k-1})+D_{k-1}-V_{k-1}$, where
$B_k$, $D_k$ and $V_k$ are nonnegative random variables, and $B_k$
and $D_k$ are summable. Thus by Theorem \ref{siegmund.thm}
%
\begin{equation}
\label{eqn17} (\overline{Y}_k-c_k)^2 \stackrel{k \rightarrow
\infty}{\longrightarrow} W \qquad\mbox{a.s.},
\end{equation}
where $ 0 \le W < \infty$ is a random variable. Thus
$|\overline{Y}_k-c_k| \rightarrow_{k \rightarrow\infty} \sqrt{W}$
a.s.

It remains to show that when $W \ne0$, $\overline{Y}_k-c_k$ cannot
jump between $\sqrt{W}$ and $-\sqrt{W}$ an infinite number of times.
It will then follow that the limit exists and is either $\sqrt{W}$
or $-\sqrt{W}$. Recall that
$\overline{Y}_k-\overline{Y}_{k-1}=\frac{Z_k}{k}$, and that by
(\ref{eqn4}) $0< c_k-c_{k-1}<\frac{\gamma}{k}$, for some $\gamma>0$.

Take expectations on both sides of the inequality in (\ref{eqn9}).
Then
\[
P\{\overline{Y}_k-\overline{Y}_{k-1}>\varepsilon\}=P\biggl\{\frac{Z_k}{k} >
\varepsilon\biggr\}=P\biggl\{\frac{Z_k^2}{k^2} > {\varepsilon^2}\biggr\} \le
\frac{C}{{\varepsilon}^2 k^{1+\delta}}.
\]
Thus by the Borel--Cantelli lemma
$P\{\overline{Y}_k-\overline{Y}_{k-1} > \varepsilon$ infinitely
often$\}=0$. This implies $P\{
|(\overline{Y}_k-c_k)-(\overline{Y}_{k-1}-c_ {k-1})| >2 \varepsilon$
infinitely often$\}=0$. Thus, if $\sqrt{W}>\varepsilon$,
$\overline{Y}_k-c_k$ cannot jump between $\sqrt{W}$ and $-\sqrt{W}$
an infinite number of times, that is, $\overline{Y}-c_k$ will converge
a.s. to $\sqrt{W}$ or $-\sqrt{W}$. Since for $W>0$, there always
exists a small enough $\varepsilon>0$ such that $W-\varepsilon>0$, it
follows that $\overline{Y}_k-c_k$ converges. Clearly on the set
where $\{W=0\}$ the statement $(\overline{Y}_k-c_k)^2 \rightarrow0$
is equivalent to $\overline{Y}_k-c_k \rightarrow0$.
\end{pf*}
\begin{pf*}{Proof of Corollary \ref{corr31}}
The expected overshoot given $X>a$ is
\[
f(a)=\frac{\int_a^\infty
e^{-x^\alpha}\,dx}{e^{-a^\alpha}}=\frac{1}{\alpha}\frac{\int_{a^\alpha
}^\infty
y^{1/\alpha-1}e^{-y}\,dy}{e^{-a^\alpha}}.
\]
The right-hand side follows by change of variables to $y=x^\alpha$.
But in Abramowitz and Stegun \cite{as}, page 263,
\[
\frac{\int_x^\infty
t^{\nu-1}e^{-t}\,dt}{e^{-x}}=x^{\nu-1} \biggl[1+\frac{\nu-1}{x}+O \biggl(\frac
{1}{x^2} \biggr) \biggr]\qquad
\mbox{as $x \rightarrow\infty$}.
\]
This implies that
\[
f(a)=\frac{1}{\alpha}a^{1-\alpha} \biggl[1+\frac{1/\alpha-1}{a^\alpha}+O
\biggl(\frac{1}{a^{2\alpha}} \biggr) \biggr]\qquad
\mbox{as $a \rightarrow\infty$}.
\]
Equation (\ref{overshoot}) implies
\[
G'(a)=\frac{1}{f(a)}=\alpha
a^{\alpha-1} \biggl[1+\frac{1-1/\alpha}{a^\alpha}+O
\biggl(\frac{1}{a^{2\alpha}} \biggr) \biggr] \qquad\mbox{as $a \rightarrow
\infty$}.
\]
Integrating both sides results in $G(a)=a^\alpha+(\alpha-1)\log
a+O(1)$ for large $a$ since the remainder term $O( a^{-(\alpha+1)})$
has finite integral.

Now for any $\varepsilon>0$, if $a$ is sufficiently large, then
$G(a^{1/\alpha}+\varepsilon)=(a^{1/\alpha}+\varepsilon)^\alpha+(\alpha-1)\log
(a^{1/\alpha}+\varepsilon)+O(1)$. But
$(a^{1/\alpha}+\varepsilon)^\alpha=a\{1+\frac{\varepsilon}{a^{1/\alpha}}\}
^\alpha=a+\alpha\varepsilon
a^{1-1/\alpha}+$ smaller order terms. Hence if $\alpha>1$ and $a$
is sufficiently large, then $G(a^{1/\alpha}+\varepsilon)>a$. A similar
argument shows that if $\alpha>1$ and $a$ is sufficiently large, then
$G(a^{1/\alpha}-\varepsilon)<a$. Therefore $\lim_{a\rightarrow\infty}
[G^{-1}(a)-a^{1/\alpha}]=0$.
\end{pf*}
\begin{pf*}{Proof of Theorem \ref{thm32}}
Upon taking expectations on both sides of (\ref{eqn7}) we obtain
%
\begin{eqnarray}\label{eqn3b1}\qquad
E[(\overline{Y}_k-c_k)^2]
&=&E[\underbrace{(\overline{Y}_{k-1}-c_{k-1})^2}_{\mathrm{(i)}}]+E \biggl(\underbrace
{\frac{E[Z^2(\overline{Y}_{k-1})|\mathcal{F}_{k-1}]}{k^2}}_{\mathrm{(ii)}} \biggr)\nonumber\\
&&{}+{}\underbrace{(c_{k-1}-c_k)^2}_{\mathrm{(iii)}}{}+{}\underbrace{\frac{2E[f(\overline
{Y}_{k-1})]}{k}(c_{k-1}-c_k)}_{\mathrm{(iv)}}\\
&&{}+{}\underbrace{2E \biggl((\overline{Y}_{k-1}-c_{k-1}) \biggl[\frac{f(\overline
{Y}_{k-1})}{k}+(c_{k-1}-c_k) \biggr] \biggr)}_{\mathrm{(v)}}.\nonumber
\end{eqnarray}
All we need to do is subtract $E[(\overline{Y}_{k-1}-c_{k-1})^2]$ on
both sides and sum. If a term remaining on the right-hand side is
positive then we need to show that it is summable. If a term is
negative it must be summable as the term on the left-hand side is
nonnegative. Hence we see that terms (ii), (iii) and (iv) cause no
trouble. The only term of concern is (v). But the expectation
(integral over the density of $\overline{Y}_k$) can be divided into
an integral over three regions: (i) $\overline{Y}_{k-1} \le c_{k-1}$,
(ii) $\overline{Y}_{k-1} \ge G^{-1}(u_k)$ and (iii) $c_{k-1}
<\overline{Y}_{k-1}<G^{-1}(u_k)$. As in the proof of Theorem
\ref{evalmost.thm}, the integrand for regions (i) and (ii) is negative
and over the third region it is positive, but can be dealt with in
the same way as in the proof of Theorem \ref{evalmost.thm} by use of
Corollary \ref{corollary2.1}. The last two statements of the theorem
follow.
\end{pf*}

\subsection{\texorpdfstring{Proofs for Section \protect\ref{sec4}}{Proofs for Section 4}}

\mbox{}

\begin{pf*}{Proof of Theorem \ref{thm41.thm}}

\textit{Proof of} (i).
%
\begin{equation}
\label{app41.eqn}\quad
E(B_k|\mathcal{F}_{k-1})=B_{k-1} \biggl[ \biggl(\frac{k-1}{k} \biggr)^{\beta-1} \biggl(1+\frac
{\beta-1}{k} \biggr) \biggr]+\frac{f(\beta\overline{Y}_{k-1})}{k^\beta}.
\end{equation}
Thus
%
\begin{equation}
\label{add41.eqn}
E(B_k|\mathcal{F}_{k-1})=B_{k-1} \biggl(1+O \biggl(\frac{1}{k^2} \biggr) \biggr)+\frac{f(\beta
(k-1)^{\beta-1}B_{k-1})}{k^\beta},
\end{equation}
where $O (\frac{1}{k^2} )>0$.
Thus
$E(B_k|\mathcal{F}_{k-1})>B_{k-1} (1+O (\frac{1}{k^2} ) )$
which implies that $B_k$ converges, to a finite or infinite limit.

Suppose first that the limit is infinite. Then there exist $k_0$ and
$D>1/\beta$ such that for all $k>k_0$, $B_{k-1}>D$. But then from
(i), for $k>k_0$
\begin{eqnarray*}
\frac{f(\beta(k-1)^{\beta-1}B_{k-1})}{k^\beta}&<&B_{k-1}\frac{c\beta
(k-1)^{\beta-1}}{k^\beta[\log(\beta
B_{k-1})+(\beta-1)\log(k-1)]^{1+\varepsilon}}\\
&<&B_{k-1}\frac{c\beta}{k[\log(\beta
D)+(\beta-1)\log(k-1)]^{1+\varepsilon}}\\
&<&B_{k-1}\frac{c\beta}{(\beta
-1)^{1+\varepsilon}
k[\log(k-1)]^{1+\varepsilon}}.
\end{eqnarray*}
But the term multiplying $B_{k-1}$ on the right is summable, which
implies that (\ref{add41.eqn}) satisfies the condition of Theorem
\ref{siegmund.thm}, and hence $B_k$ converges to a finite limit.
This contradiction implies that $B_k$ converges to a finite r.v.
a.s.

\textit{Proof of} (ii).\quad
Now suppose that $B_k$ converges a.s. and $\lim EB_k<\infty$. It
follows from (\ref{eqn41}) and as in (\ref{app41.eqn}) that $B_k$
can be written as
\[
B_k=B_{k-1} \biggl[1+O \biggl(\frac{1}{k^2} \biggr) \biggr]+\frac{Z(\beta\overline
{Y}_{k-1})}{k^\beta},
\]
where $O(1/k^2)$ is positive. It follows that $B_k>B_{k-1}$, so that
the limit is positive. Since the support of the observations is not
bounded, in a different realization one could obtain a higher value.
Hence the limit is a nondegenerate positive random variable. Set
$B_0=0$. Then
%
\begin{equation}
\label{eqn410}
B_k=\sum_{j=1}^k(B_j-B_{j-1})=O(1)\sum_{j=1}^k\frac{B_{j-1}}{j^2}+\sum
_{j=2}^k\frac{Z(\beta\overline{Y}_{j-1})}{j^\beta}.
\end{equation}

Since $B_k$ converges a.s. $\lim B_k$ exists and is finite a.s.
Taking expectations and limits as $k \rightarrow\infty$ on both
sides of (\ref{eqn410}) and noting that $EB_k$ is assumed to be
bounded implies that
$\sum_{j=1}^\infty\frac{Z(\beta\overline{Y}_{j-1})}{j^\beta} <
\infty$. This in turn implies
$\sum_{j=1}^\infty\frac{f(\beta\overline{Y}_{j-1})}{j^\beta}<
\infty$.
%
Since $B_k$ converges a.s. to a random variable $W_\beta$ for
$0<\varepsilon<W_\beta$ and a (random) $k_0$, we have for all
$k-1>k_0$,
$(W_\beta-\varepsilon)(k-1)^{\beta-1}<\overline{Y}_{k-1}<(W_\beta+\varepsilon
)(k-1)^{\beta-1}$.
If $f$ is increasing
%
\begin{eqnarray}
\label{eqn49}
\infty&>&\sum_{k=k_0}^\infty\frac{f(\beta
\overline{Y}_{k-1})}{k^\beta}
>\sum_{k=k_0}^\infty\frac{f(\beta(W_\beta-\varepsilon)(k-1)^{\beta
-1})}{k^\beta}\nonumber\\[-8pt]\\[-8pt]
&>& \biggl(\frac{1}{2} \biggr)^\beta\sum_{k=k_0}^\infty\frac{f(A(k-1)^{\beta
-1})}{(k-1)^\beta},\nonumber
\end{eqnarray}
where the inequality follows since $ (\frac{k-1}{k} )^\beta
> (\frac{1}{2} )^\beta$, and where we have written
$A=\beta(W_\beta-\varepsilon)$.

Finally,
\begin{eqnarray*}
\sum_{k=k_0}^\infty\frac{f(A(k-1)^{\beta-1})}{(k-1)^\beta}
&>&\int_{k_0-2}^\infty\frac{1}{(x+1)^\beta}f(Ax^{\beta-1})\,dx\\
&>& \biggl(\frac{k_0-2}{k_0-1} \biggr)^\beta
\int_{k_0-2}^\infty\frac{f(Ax^{\beta-1})}{x^\beta}\,dx.
\end{eqnarray*}
By change of variable to $y=Ax^{\beta-1}$ the integral on the
right-hand side becomes
$\frac{A}{\beta-1}\int_{A(k_0-2)^{\beta-1}}^\infty\frac{f(y)}{y^2}\,dy$.
This integral is therefore finite by (\ref{eqn49}).
%
\end{pf*}
%
%
\begin{pf*}{Proof of Proposition \ref{prop41}}
In a manner similar to the end of the proof in Theorem
\ref{thm41.thm}, it can be shown that if
${\int_C^\infty}\frac{\Psi(y)}{y^2}\,dy$ diverges, then ${\lim_{n
\rightarrow\infty}\sum_{k=k_0}^n}\frac{1}{k^\beta} 
\Psi(\gamma
k^{\beta-1})$ also diverges.

Note that
%
\begin{equation}
\label{ext1.eqn} B_k=B_{k-1} \biggl[1+
O \biggl(\frac{1}{k^2} \biggr) \biggr]+\frac{Z^*(\beta\overline{Y}_{k-1})}{k^\beta}\Psi
(\beta\overline{Y}_{k-1}).
\end{equation}
Let $F_k$ be the c.d.f. of $Z^*(\beta\overline{Y}_{k-1})$ conditional
on $\overline{Y}_{k-1}$. Let $F_V$ be the c.d.f. of $V$. Let
$U_1,U_2,\ldots\sim U[0,1]$ i.i.d. Define $V_k=F_V^{-1}(U_k)$ (so
$V_i$ are i.i.d. with c.d.f. $F_V$). Clearly, $V_k \le F_k^{-1}(U_k)$
conditional on $\overline{Y}_{k-1}$ once $\beta\overline{Y}_{k-1}
\ge a_0$ (which will happen with probability 1).

It follows that one can imbed the sequence $Y_1,Y_2,\ldots$ in a
probability space where $V_1,V_2,\ldots$ are i.i.d. with c.d.f. $F_V$
and
\[
V_i \le Z^*(\beta\overline{Y}_{k-1}) \qquad\mbox{for all $i$}\qquad\mbox{such that
$\beta\overline{Y}_{i-1} \ge a_0$}
\]
conditional on $\overline{Y}_{k-1}$. Define $V_i^*=c1(V_i>c)$ for
some $c$ such that $P(V_i>c)>0$. Clearly, $V_i^* \le
Z^*(\beta\overline{Y}_{k-1})$. Note that $V_i^*$ is $c$ times a
Bernoulli random variable. Now
%
\begin{equation}
\label{ext2.eqn}
\frac{Z^*(\beta\overline{Y}_{k-1})}{k^\beta}\Psi(\beta\overline{Y}_{k-1})
\ge V_k^*\frac{\Psi(\beta\overline{Y}_{k-1})}{k^\beta}.
\end{equation}
Recall that
%
\[
\overline{Y}_k
\ge\frac{\beta-1+k}{k}\overline{Y}_{k-1}
\]
so that for a constant $a_1$ that is independent of $Y_1$ and $k$,
%
\begin{equation}
\label{ext3.eqn} \overline{Y}_k >
Y_1\prod_{j=2}^k\frac{\beta-1+j}{j}
\ge Y_1a_1k^{\beta-1}.
\end{equation}
Hence (for $k$ such
that $\beta\overline{Y}_{k-1} \ge a_0$ and $a_1$ a constant) since
$\Psi(a)$ increases in $a$
%
\begin{equation}
\label{ext4.eqn}
\frac{Z^*(\beta\overline{Y}_{k-1})}{k^\beta}\Psi(\beta\overline{Y}_{k-1})
\ge V_k^*\frac{\Psi(\beta Y_1a_1(k-1)^{\beta-1})}{k^\beta}.
\end{equation}
Finally, condition on $Y_1$ and denote
\[
c_k=\frac{\Psi(\beta Y_1a_1(k-1)^{\beta-1})}{k^\beta}.
\]
By (\ref{prop41.eqn}) and what we showed above, $\lim_{n \rightarrow
\infty}\sum_{k=1}^nc_k =\infty$.
It is a straightforward application of Kolmogorov's three-series
theorem (cf. Feller \cite{feller}, page 317) that
%
\begin{equation}
\label{ext5.eqn} \lim_{n \rightarrow\infty} \sum_{k=1}^nV^*_kc_k
=\infty\qquad\mbox{a.s.}
\end{equation}
Putting
(\ref{ext4.eqn}) and (\ref{ext5.eqn})
together obtains that
$\frac{Z^*(\beta\overline{Y}_{k-1})}{k^\beta}\Psi(\beta\overline{Y}_{k-1})$
is not summable. This and (\ref{ext1.eqn}) imply that $\lim_{k
\rightarrow\infty}B_k=\infty$ a.s.
\end{pf*}
%
%
\begin{pf*}{Proof of Theorem \ref{lemm42}}
$EB_k$ converges to a finite limit by (\ref{eqn48n}) since $f$ is
bounded by assumption (b)
%
\begin{eqnarray}
\label{eqnew1}\qquad \operatorname{Var} B_k&=&\operatorname{Var} \biggl(\frac
{k-1+\beta}{{k^\beta}}\overline{Y}_{k-1}+\frac{Z(\beta
\overline{Y}_{k-1})}{k^\beta} \biggr)\nonumber\\
&=& \biggl[\frac{(k-1+\beta)^2(k-1)^{2(\beta-1})}{k^{2\beta}} \biggr]
\operatorname{Var}B_{k-1}+\frac{\operatorname{Var}(Z(\beta\overline{Y}_{k-1}))}{k^{2\beta}}\\
&&{}+
2\frac{k-1+\beta}{k^{2\beta}}\operatorname{Cov}(\overline{Y}_{k-1},Z(\beta\overline
{Y}_{k-1}) ).\nonumber
\end{eqnarray}
We shall treat each of the three terms in (\ref{eqnew1})
separately.

\begin{longlist}
\item It is easily\vspace*{2pt} seen (by taking $\log$) that
the value in the square bracket is $1+O (\frac{1}{k^2} )$.\vspace*{2pt}

\item From condition (c) and the convergence of $EB_k$ to a finite
limit
\[
\frac{\operatorname{Var}(Z(\beta\overline{Y}_{k-1}))}{k^{2\beta}}<\frac{EZ^2(\beta
\overline{Y}_{k-1})}{k^{2\beta}}
<\frac{c\beta E\overline{Y}_{k-1}}{k^{2\beta}} <\frac{c\beta(\lim
EB_k+\varepsilon)}{k^{\beta+1}}.
\]
Thus the second term in the right-hand side of (\ref{eqnew1}) is
summable.

\item We now show that the third term in the right-hand side
of (\ref{eqnew1}) is negative or 0:
\begin{eqnarray*}
&&\operatorname{Cov}(\overline{Y}_{k-1},Z(\beta\overline{Y}_{k-1}))
\\
&&\qquad=E(\overline{Y}_{k-1}Z(\beta\overline{Y}_{k-1}))
-E(\overline{Y}_{k-1})E(Z(\beta\overline{Y}_{k-1}))\\
&&\qquad=E[\overline{Y}_{k-1}E(Z(\beta\overline{Y}_{k-1})|\mathcal
{F}_{k-1})]-E(\overline{Y}_{k-1})E[E(Z(\beta\overline{Y}_{k-1})|\mathcal
{F}_{k-1})]\\
&&\qquad=\frac{1}{\beta}E[\beta\overline{Y}_{k-1}f(\beta\overline
{Y}_{k-1})]-\frac{1}{\beta}E(\beta\overline{Y}_{k-1})Ef(\beta\overline
{Y}_{k-1})\\
&&\qquad=\frac{1}{\beta}\operatorname{Cov}(\beta\overline{Y}_{k-1},f(\beta\overline{Y}_{k-1}))
\le0,
\end{eqnarray*}
where the last inequality follows from (b). It follows that
(\ref{eqnew1}) satisfies the condition in Corollary
\ref{corollary2.1} with $z_n=\operatorname{Var} B_n$, and the result
follows.\qed
\end{longlist}
\noqed\end{pf*}

\subsection{\texorpdfstring{Proofs for Section \protect\ref{sec5}}{Proofs for Section 5}}

\mbox{}

\begin{pf*}{Proof of Theorem \ref{thm51}}
Let $P_j=1-F(\overline{Y}_{j-1})$. We shall use Theorem
\ref{feller.thm} conditionally on the sequence $\{\overline{Y}_k\}$.
Let $b_j=\sum_{i=1}^j P_i^{-1}$ and $Q_i=T_i-T_{i-1}-P_i^{-1}$ with
$T_0\equiv0$. Obviously, the sequence $\{b_j\}_{j=1}^\infty$
satisfies the first condition of Theorem \ref{feller.thm}.

Note that conditional on the sequence $\{P_j\}$ the distribution of
$T_i-T_{i-1}$ is Geometric $(P_i)$ and these differences are
conditionally independent of each other. Hence
$\{Q_n\}_{n=1}^\infty$ is a sequence of conditionally independent
random variables with zero expectation and variance $(1-P_n)/P_n^2$.
We shall show that the second condition of Theorem \ref{feller.thm}
holds
\[
\sum_{n=1}^\infty
E(Q_n^2/b_n^2)=\sum_{n=1}^\infty\frac{1-P_n}{P_n^2}\bigg/\Biggl(\sum_{j=1}^nP_j^{-1}\Biggr)^2
< \sum_{n=1}^\infty\frac{1}{P_n^2}\bigg/\Biggl(\sum_{j=1}^n P_j^{-1}\Biggr)^2.
\]

It therefore suffices to show that for all $n\ge n_0$
%
\begin{equation}
\label{eqn53} \sum_{j=0}^n\frac{P_{n+1}}{P_{j+1}} \ge An^{1/2}\log n
\end{equation}
for some $A>0$. We shall actually show that for any $0<\varepsilon<
1/4$ there exists $j_0$ such that for all $n \ge j \ge j_0$
%
\begin{equation}
\label{eqn54} \frac{P_{n+1}}{P_{j+1}} >
\frac{j^{1-\varepsilon}}{n^{1+\varepsilon}}.
\end{equation}
From (\ref{eqn54}) it is immediate that (\ref{eqn53}) holds, since
\[
\sum_{j=0}^n\frac{P_{n+1}}{P_{j+1}}>\sum_{j=j_0}^n\frac{P_{n+1}}{P_{j+1}}
\ge\frac{1}{n^{1+\varepsilon}}\sum_{j=j_0}^nj^{1-\varepsilon} \ge
D\frac{n^{2-\varepsilon}-j_0^{2-\varepsilon}}{n^{1+\varepsilon}}>An^{1-2\varepsilon}.
\]
Note that for $H(x)=x^\alpha+h(x)$ for $\alpha>0$ and $h$ that
satisfies (\ref{eqn1n}), we have by Theorem \ref{thm3n},
$\overline{Y}_j=(\log j)^{1/\alpha}(1+\varepsilon_j)$ with $\varepsilon_j
\stackrel{j \rightarrow\infty}{\longrightarrow} 0$. Thus
\[
H(\overline{Y}_j)=(\log
j)(1+\varepsilon_j)^\alpha\biggl[1+\frac{h((\log
j)^{1/\alpha}(1+\varepsilon_j))}{(\log j)(1+\varepsilon_j)^\alpha} \biggr]
\]
and since $h(x)/x^\alpha\stackrel{x \rightarrow
\infty}{\longrightarrow} 0$ it follows that for any $\varepsilon>0$
there exists $j_0$ such that for all $j>j_0$
\[
(1+\varepsilon)\log j >H(\overline{Y}_j)>(1-\varepsilon)\log j,
\]
which implies, since $[1-F(\overline{Y}_j)]^{-1}=\exp
H(\overline{Y}_j)$, that
%
\begin{equation}
\label{eqnn71}
j^{1+\varepsilon}>[1-F(\overline{Y}_j)]^{-1}>j^{1-\varepsilon}.
\end{equation}
Thus (\ref{eqn54}) follows.

Note that here $S_n$ of Theorem \ref{feller.thm} equals
$T_n-\sum_{i=1}^nP_i^{-1}$, thus $b_n^{-1}S_n \rightarrow0$ a.s. is
equivalent to $T_n^*-1 \rightarrow0$ a.s. Since this result holds
for any conditioning sequence~$\{\overline{Y}_k\}$, it holds
unconditionally.
\end{pf*}
%
%
\begin{pf*}{Proof of Corollary \ref{cor51}}
In (\ref{eqnn71}) take any $\varepsilon>0$. Hence for some positive
constants $c_1, c_2, c^*_1, c^*_2$ and all $k$ large enough
\[
c^*_2k^{2+\varepsilon}>c_2\sum_{j=1}^kj^{1+\varepsilon}>\sum
_{j=1}^k[1-F(\overline{Y}_j)]^{-1}>c_1\sum_{j=1}^kj^{1-\varepsilon
}>c^*_1k^{2-\varepsilon}.
\]
Since $\frac{T_k}{\sum_{j=1}^k[1-F(\overline{Y}_j)]^{-1}}
\rightarrow1$ a.s. for $k$ large enough and $c_1^{**}$ a positive
constant
\begin{eqnarray*}
\frac{T_k}{k^{2-\delta}} & = &
\frac{T_k}{\sum_{j=1}^k[1-F(\overline{Y}_j)]^{-1}}
\frac{\sum_{j=1}^k[1-F(\overline{Y}_j)]^{-1}}{k^{2-\delta}}\\
& > &
c_1^{**}k^{\delta-\varepsilon} \rightarrow\infty\qquad\mbox{a.s.
if }
\delta>\varepsilon.
\end{eqnarray*}
The proof for $\frac{T_k}{k^{2+\delta}}$ follows in a similar
manner.
\end{pf*}
\begin{pf*}{Proof of Theorem \ref{betasmall}}

\textsc{Proof of (i)}.\quad We shall (again) use Theorem \ref{feller.thm} and
show (\ref{eqn53})
where $P_j=1-F(\beta\overline{Y}_{j-1})$.
Assume that\vspace*{1pt} $k_0$ (random) is such that for all $k \ge k_0$, $Z_k
<\gamma\overline{Y}_{k-1}$. Such a $k_0$ exists with probability
one by Lemma \ref{lem6}. Then for $k>k_0$,
\[
\overline{Y}_k=\overline{Y}_{k-1}+\frac{Z_k+(\beta-1)\overline{Y}_{k-1}}{k}
\le\overline{Y}_{k-1} \biggl(1+\frac{\gamma+\beta-1}{k} \biggr).
\]
Thus
%
\begin{equation}
\label{eqn64}
\overline{Y}{}^\alpha_k\le
\overline{Y}{}^\alpha_{k-1}\biggl(1+\frac{\gamma+\beta-1}{k} \biggr)^\alpha
\le\overline{Y}{}^\alpha_{k-1}\biggl(1+\frac{d}{k} \biggr),
\end{equation}
where $d=(\gamma+\beta-1)\rho_{u\alpha}$ and $\rho_{u\alpha}$ (and
for later purposes $\rho_{l\alpha}$) is defined by
%
\begin{equation}
\label{rho} 1+\rho_{l\alpha}x \le(1+x)^\alpha\le1+\rho_{u\alpha}x\qquad
\mbox{for $0 \le x \le1$},
\end{equation}
where
$\rho_{l\alpha}=\alpha$ and $\rho_{u\alpha}=2^\alpha-1$ when $\alpha
\ge1$, while $\rho_{l\alpha}=2^\alpha-1$ and
$\rho_{u\alpha}=\alpha$ when $\alpha<1$. We have used the inequality
$(1+x)^\alpha\le1+(2^\alpha-1)x$, valid for all $\alpha\ge1$ and
$0 \le x \le1$.
We can therefore write, using (\ref{eqn64}),
%
\begin{equation}
\label{eqn67}
\frac{P_{k+1}}{P_k}= \exp
\{-\beta^\alpha(\overline{Y}{}^\alpha_k-\overline{Y}{}^\alpha_{k-1})\}
\ge
\exp\biggl\{-\beta^\alpha\frac{d}{k}\overline{Y}{}^\alpha_{k-1}\biggr\}.
\end{equation}
Now let $k_1 \ge k_0$ be so large that for all $k \ge k_1$,
$\overline{Y}_k<(W+\varepsilon)k^{\beta-1}$, which exists, by Theorem
\ref{meanstretch}. Then we can continue the inequality in
(\ref{eqn67}), by
\[
\frac{P_{k+1}}{P_k} >
\exp\biggl\{-\beta^\alpha\frac{d}{k}(W+\varepsilon)^\alpha
k^{\alpha(\beta-1)} \biggr\}
=\exp\bigl\{-Bk^{\alpha(\beta-1)-1} \bigr\}.
\]
To simplify notation let
%
\begin{equation}
\label{eqn68} \tau= \alpha(\beta-1)-1,
\end{equation}
thus $\tau> -1$. For $j >k>k_1$, we have
\[
\frac{P_{n+1}}{P_{j+1}}=\prod_{k=j+1}^n\frac{P_{k+1}}{P_k}
>\exp\Biggl\{-B\sum_{k=j+1}^n k^\tau
\Biggr\}>\exp\biggl\{-\frac{B}{\tau+1} \bigl(n^{\tau+1}-(j+1)^{\tau+1}
\bigr) \biggr\}.
\]
Thus
\[
\sum_{j=1}^n\frac{P_{n+1}}{P_{j+1}}
>\sum_{j=k_1}^n\frac{P_{n+1}}{P_{j+1}}
>e^{-[{B}/({\tau+1})]n^{\tau+1}}
\sum_{j=k_1}^ne^{[{B}/({\tau+1})](j+1)^{\tau+1}}.
\]
But
\[
\sum_{j=k_1}^ne^{[{B}/({\tau+1})](j+1)^{\tau+1}}>\int
_{k_1+1}^{n}e^{[{B}/({\tau+1})]x^{\tau+1}}\,dx,
\]
thus
%
\begin{equation}
\label{eqn69}
\sum_{j=1}^n\frac{P_{n+1}}{P_{j+1}}>\int_{k_1+1}^{n}e^{[{B}/({\tau
+1})]x^{\tau+1}}\,dx/e^{[{B}/({\tau+1})]n^{\tau+1}}.
\end{equation}
We would like the right-hand side of (\ref{eqn69}), divided by
$n^{1/2 +\varepsilon}$ for some (small) $\varepsilon>0$, to tend to a
nonzero limit in order for (\ref{eqn53}) to hold. Thus consider, by
use of l'H\^{o}pital's rule, the limit as $y \rightarrow\infty$
of
\begin{eqnarray*}
q(y)&=&\frac{\int_{k_1+1}^ye^{Ax^{\tau+1}}\,dx}{y^\delta
e^{Ay^{\tau+1}}}\qquad \mbox{where $A>0$ is any constant},
\\
\lim_{y \rightarrow\infty} q(y)&=&\lim_{y \rightarrow
\infty}\frac{e^{Ay^{\tau+1}}}{e^{Ay^{\tau+1}}(\delta
y^{\delta-1}+A(\tau+1)y^{\tau+\delta})},
\end{eqnarray*}
which is finite when $\tau+\delta=0$ and tends to $\infty$ when
$\tau+\delta<0$. Now for $\delta=1/2$, by (\ref{eqn68}) we get a
finite limit when $\alpha(\beta-1)-1+1/2=0$, that is,
$\beta=1+1/(2\alpha)$. Thus for $\beta<1+1/(2\alpha)$ there will
exist an $\varepsilon>0$ such that the value of
$\sum_{j=0}^{n+1}\frac{P_{n+1}}{P_{j+1}}>n^{1/2+\varepsilon}$, and the
result (i) follows.\vspace*{1pt}

\textsc{Proof of} (ii).\quad Let
%
\begin{equation}
\label{eqnm62}
\gamma(k)=\frac{1}{\sum_{j=1}^{k-1}e^{(\beta\overline{Y}_j)^\alpha}}
\end{equation}
in (\ref{eqnm61}). Then clearly $\gamma(k) \rightarrow0$. We shall
show later that
$ [(1+o_k(1))t\gamma(k)\times\break e^{(\beta\overline{Y}_{j-1})^\alpha} ]$
of (\ref{eqnm61}) is arbitrarily close to $0$ for $2 \le j \le k$
for all sufficiently large $k$ and $\beta<1+1/\alpha$. It suffices
to show this for $j=k$. We can then write, using (\ref{eqnm61}) and
(\ref{eqnm62}),
\begin{eqnarray*}
-(1+\varepsilon
)t&=&-(1+\varepsilon)\sum_{j=2}^kt\gamma(k)e^{(\beta\overline
{Y}_{j-1})^\alpha}
< \log Ee^{-t\gamma(k)\tilde{T}_k}\\
&<&-(1-\varepsilon)\sum_{j=2}^kt\gamma(k)e^{(\beta\overline
{Y}_{j-1})^\alpha}=-(1-\varepsilon
)t.
\end{eqnarray*}
It follows that $\lim_{k \rightarrow
\infty}E (e^{-t\gamma(k)T_k} )=e^{-t}$, which is the
desired result. We still must show that
$ [(1+o_k(1))t\gamma(k)e^{(\beta\overline{Y}_{j-1})^\alpha} ]$
of (\ref{eqnm61}) is arbitrarily close to $0$ for $j= k$ for all
sufficiently large $k$ and $\beta<1+1/\alpha$. Let $\rho_{l\alpha}$
be defined by (\ref{rho})
\[
\overline{Y}{}^\alpha_j-\overline{Y}{}^\alpha_{j-1}>\frac{\rho_{l\alpha
}(\beta-1)\overline{Y}{}^\alpha_{j-1}}{j}>\frac{\rho_{l\alpha}(\beta
-1)j^{(\beta-1)\alpha}W^\alpha(1-\varepsilon)}{j}
\]
for all $j$ sufficiently large, where by Theorem \ref{meanstretch}
$\lim\frac{\overline{Y}_{j-1}}{(j-1)^{\beta-1}}=W>0$. Thus
\begin{eqnarray*}
\gamma(k)e^{(\beta\overline{Y}_{k-1})^\alpha}&=&\frac{1}{\sum
_{j=1}^{k-1}e^{-\beta^\alpha(\overline{Y}{}^\alpha_{k-1}-\overline
{Y}{}^\alpha_j)}}=\frac{1}{\sum_{j=1}^{k-1}e^{-\beta^\alpha\sum
_{i=j+1}^{k-1}(\overline{Y}{}^\alpha_i-\overline{Y}{}^\alpha_{i-1})}}\\
&<&\frac{1}{\sum_{j=j_0}^{k-1}e^{-D\sum_{i=j+1}^{k-1}i^{(\beta-1)\alpha
-1}}}\rightarrow
0,
\end{eqnarray*}
for suitable large $j_0$, as long as $(\beta-1)\alpha-1<0$, that is,
$\beta<1+1/\alpha$ [where we have let
$D=\rho_{l\alpha}\beta^\alpha(\beta-1)W^\alpha(1-\varepsilon)$].

\textsc{Proof of} (iii).\quad Here let $\gamma(k)=e^{-\beta^\alpha\overline
{Y}{}^\alpha_{k-1}}$.
With this $\gamma(k)$, (\ref{eqnm61}) becomes
%
\begin{eqnarray}
\label{eqnb64}\quad
\log
Ee^{-t\gamma(k)\tilde{T}_k}&=&
-\sum_{j=2}^k\log\bigl[1+\bigl(1+o_k(1)\bigr)te^{-\beta^\alpha(\overline
{Y}{}^\alpha_{k-1}-\overline{Y}{}^\alpha_{j-1})} \bigr]\nonumber\\
&=&-\log\bigl[1+\bigl(1+o_k(1)\bigr)t\bigr]\\
&&{}-\sum_{j=2}^{k-1}\log\bigl[1+\bigl(1+o_k(1)\bigr)te^{-\beta^\alpha\sum
_{i=j}^{k-1}(\overline{Y}{}^\alpha_i-\overline{Y}{}^\alpha_{i-1})}
\bigr].\nonumber
\end{eqnarray}
Now for some $D>0$ (dependent on $\{\overline{Y}_k\})$
\[
0<e^{-\beta^\alpha(\overline{Y}{}^\alpha_i-\overline{Y}{}^\alpha_{i-1}
)}<e^{\frac{-\beta^\alpha\rho_{l\alpha}(\beta-1)\overline
{Y}{}^\alpha_{i-1}}{i}}<e^{-Di^{\alpha(\beta-1)-1}}.
\]
Thus
\[
e^{-\beta^\alpha(\overline{Y}{}^\alpha_{k-1}-\overline{Y}{}^\alpha_{j-1}
)}<e^{-D\int_j^{k-1}x^\nu
\,dx}=e^{-{D}/({\nu+1})[(k-1)^{\nu+1}-j^{\nu+1}]},
\]
where $\nu=\alpha(\beta-1)-1>0$, that is, $\beta>1+1/\alpha$.
But
\begin{eqnarray*}
\lim_{k \rightarrow
\infty}\sum_{j=2}^{k-1}e^{-{D}/({\nu+1})[(k-1)^{\nu+1}-j^{\nu
+1}]}&=&\lim_{k
\rightarrow\infty}
\frac{\int_2^ke^{Dx^{\nu+1}}\,dx}{e^{Dk^{\nu+1}}}\\
&\stackrel{\mathrm{l}\mbox{\fontsize{8.36}{10}\selectfont{'
H\^{o}pital}}}{=}&\lim_{k \rightarrow\infty}
\frac{1}{D(\nu+1)k^\nu}= 0.
\end{eqnarray*}
Since the sum in the right-hand side of (\ref{eqnb64})
tends to $0$ and
\[
\lim_{k \rightarrow
\infty}E \bigl(e^{-t\gamma(k)T_k} \bigr) \rightarrow1/(1+t),
\]
which is $Ee^{-tQ}$ where $Q \sim \operatorname{Exp}(1)$,
$\frac{T_k}{e^{\beta\alpha}\overline{Y}{}^\alpha_{k-1}}$ tends in
distribution to an exponential distribution. The above proof shows that
$\sum_{j=1}^{k-1}e^{(\beta\overline{Y}_j)^\alpha}/e^{(\beta\overline
{Y}{}^\alpha_{k-1})} \stackrel{\mathrm{a.s.}}{\longrightarrow} 1$, thus also
$T^*_k \stackrel{D}{\longrightarrow} \operatorname{Exp}(1)$.
\end{pf*}
\begin{lemma}
\label{lem6}
Let $X \sim F$ where $F$ is $\mathcal{G}_\alpha$ with $\alpha> 0$.
Let $\beta>1$ and let $Z_k$ be the random ``overshoot'' over $\beta
\overline{Y}_{k-1}$. For any $\gamma>0$, and any $0 \le\delta<
(\beta-1)\alpha$, $P(Z_k>\gamma\overline{Y}_{k-1}/k^\delta
\mbox{ infinitely often})=0$.
\end{lemma}
\begin{pf}
Consider the event $A=\{\overline{Y}_k/k^{\beta-1} \rightarrow W,
0<W<\infty\}$. We know by Theorem \ref{meanstretch} that $P(A)=1$,
and hence we shall assume that $A$ occurs. Let $A_k=\{Z_k
>\gamma\overline{Y}_{k-1}/k^\delta\}$. We shall show that
$\sum_{k=1}^\infty P(A_k)< \infty$ so that the result will follow
from the Borel--Cantelli lemma.
\begin{eqnarray*}
P(A_k|\overline{Y}_{k-1})&=&\exp\bigl\{-[(\gamma/k^\delta+\beta
)^\alpha-\beta^\alpha]\overline{Y}{}^\alpha_{k-1}\\
&&\hspace*{18.4pt}{}-\bigl[h\bigl((\gamma/k^\delta+\beta)\overline{Y}_{k-1}\bigr)-h(\beta\overline
{Y}_{k-1})\bigr]\bigr\}\\
&=&\exp
\biggl\{-\beta^\alpha\biggl( \biggl(1+\frac{\gamma/\beta}{k^\delta} \biggr)^\alpha-1 \biggr)\overline
{Y}{}^\alpha_{k-1}
-h'(Q_k)\frac{\gamma}{k^\delta}\overline{Y}_{k-1}\biggr\},
\end{eqnarray*}
where $\beta\overline{Y}_{k-1} \le Q_k \le
(\beta+\frac{\gamma}{k^\delta})\overline{Y}_{k-1}$. Write, by
(\ref{eqn5n}),
\[
\biggl|-h'(Q_k)\frac{\gamma}{k^\delta}\overline{Y}_{k-1}\biggr|=\biggl|-\frac
{h'(Q_k)}{Q_k^{\alpha-1}}\frac{\gamma}{k^\delta}\overline
{Y}_{k-1}Q_k^{\alpha-1}\biggr|=\biggl|o_k\overline{Y}{}^\alpha_{k-1}\frac{1}{k^\delta}\biggr|,
\]
where $|o_k| <\varepsilon$ for $k \ge k_0$ with large enough $k_0$
and $\varepsilon>0$ arbitrary. Note that by (\ref{ext3.eqn}) $k_0$ can be
chosen to depend on $Y_1$ only. For
$\varepsilon$ small enough this implies
\begin{eqnarray*}
P(A_k|\overline{Y}_{k-1}) & \le&
\exp\biggl\{-\beta^\alpha\biggl(\frac{\alpha\gamma/\beta}{k^\delta}
+o \biggl(\frac{1}{k^\delta} \biggr) \biggr)\overline{Y}{}^\alpha_{k-1}+|o_k|\frac
{1}{k^\delta}\overline{Y}{}^\alpha_{k-1}\biggr\}\\
& \le&
\exp\biggl\{-[\beta^{\alpha-1}\alpha\gamma-2\varepsilon]\frac{1}{k^\delta
}\overline{Y}{}^\alpha_{k-1}\biggr\}\\
&\le
&\exp\biggl\{-\frac{1}{2}\beta^{\alpha-1}\alpha\gamma\frac{1}{k^\delta
}cY_1k^{\beta-1}\biggr\}
\le\exp\{-dY_1k^{\beta-1-\delta}\}.
\end{eqnarray*}
The next to last inequality follows from (\ref{ext3.eqn}). Hence
$P(A_k|Y_1)\le\exp\{-dY_1\times\break k^{\beta-1-\delta}\}$ for $k \ge
k_0=k_0(Y_1)$. If $\delta<\beta-1$,
$\sum_{k=1}^{\infty}P(A_k|Y_1)<\infty$, so, by the Borel--Cantelli
lemma, conditional on $Y_1$, $P(A_k$ i.o.$|Y_1)=0$. But this is true
for all $Y_1$. Hence $P(A_k$ i.o.$)=0$.
\end{pf}
\begin{pf*}{Proof of Theorem \ref{thmeq}}
Let $\gamma(k)=1/\sum_{j=1}^{k-1}e^{d\overline{Y}{}^\alpha_j}$ where
$d=\beta^\alpha$. We write (\ref{eqnm61}) as
%
\begin{equation}
\label{eqn71} \log
E \bigl(e^{-t\gamma(k)\tilde{T}_k} \bigr)=-\sum_{j=1}^{k-1}\log
\bigl[1+t\bigl(1+o_k(1)\bigr)\gamma(k)e^{d\overline{Y}{}^\alpha_{k-j}} \bigr].
\end{equation}
%

Since when $R$ is exponentially distributed with mean $\mu$, $\log
E (e^{-tR} )=-{\log}(1+\mu t)$, it is sufficient to show that
the right-hand side of (\ref{eqn71}) converges, as $k \rightarrow
\infty$, to $-\sum_{j=1}^\infty\log[1+t\mu_j]$.

First consider
%
\begin{equation}
\label{eqn72}
\gamma(k)e^{d\overline{Y}{}^\alpha_{k-j}}=\frac{1}{S_{j,k}+T_{j,k}},
\end{equation}
where
$S_{j,k}=\sum_{i=k-j}^{k-1}e^{d(\overline{Y}{}^\alpha_i-\overline
{Y}{}^\alpha_{k-j})}$
and
$T_{j,k}=\sum_{i=1}^{k-j-1}e^{-d(\overline{Y}{}^\alpha_{k-j}-\overline
{Y}{}^\alpha_i)}$.
Note that $Y_i=Z_i+\beta\overline{Y}_{i-1}$ where $Z_i$ is the
amount above $\beta\overline{Y}_{i-1}$ for the $i$th item that
is kept. Hence,
\begin{eqnarray*}
\overline{Y}_i&=&\frac{(i-1)\overline{Y}_{i-1}+Z_i+\beta\overline
{Y}_{i-1}}{i}\\
&=&\overline{Y}_{i-1}+\frac{Z_i+\overline{Y}_{i-1}/\alpha}{i}\qquad
\mbox{because } \beta-1=1/\alpha.
\end{eqnarray*}
By Lemma \ref{lem6}, for all $i$ sufficiently large
\begin{eqnarray*}
\overline{Y}{}^\alpha_i&=&\overline{Y}{}^\alpha_{i-1}\biggl(1+\frac
{Z_i}{i\overline{Y}_{i-1}}+\frac{1}{\alpha i} \biggr)^\alpha\\
&=&\overline{Y}{}^\alpha_{i-1}\biggl(1+\frac{1}{i}+ \mbox{smaller
order terms} \biggr).
\end{eqnarray*}
Let $w=\lim_{ k \rightarrow\infty}\overline{Y}_k/k^{1/\alpha}$.
Therefore, $\lim_{i \rightarrow
\infty}\overline{Y}{}^\alpha_i-\overline{Y}{}^\alpha_{i-1}=w^\alpha$ and
for fixed $b$, $\lim_{i \rightarrow\infty}
\overline{Y}{}^\alpha_{i+b}-\overline{Y}{}^\alpha_i=bw^\alpha$. This
implies
%
\begin{eqnarray}
\label{eqn73} \lim_{k \rightarrow\infty} S_{j,k}&=&\lim_{k
\rightarrow
\infty}\sum_{i=k-j}^{k-1}e^{d(\overline{Y}{}^\alpha_i-\overline
{Y}{}^\alpha_{k-j})}
=\lim_{k \rightarrow
\infty}\sum_{l=0}^{j-1}e^{d(\overline{Y}{}^\alpha_{k-j+l}-\overline
{Y}{}^\alpha_{k-j})}\nonumber\\[-8pt]\\[-8pt]
&=&\sum_{l=0}^{j-1}e^{dlw^\alpha} =\frac{e^{dw^\alpha
j}-1}{e^{dw^\alpha}-1}.\nonumber
\end{eqnarray}

For any $\varepsilon>0$ there exists $m$ such that
$(1-\varepsilon)w^\alpha\le\frac{\overline{Y}{}^\alpha_i}{i+1} \le
(1+\varepsilon)w^\alpha$ for all $i \ge m$. This implies
\[
\lim_{k \rightarrow\infty}T_{j,k}=\lim_{k \rightarrow
\infty}\sum_{i=1}^{m-1}e^{-d(\overline{Y}{}^\alpha_{k-j}-\overline
{Y}{}^\alpha_i)}+\lim_{k
\rightarrow
\infty}\sum_{i=m}^{k-j-1}e^{-d(\overline{Y}{}^\alpha_{k-j}-\overline
{Y}{}^\alpha_i)}.
\]
Fix $m$. Then the first limit on the right-hand side is clearly zero
since $\overline{Y}_{k-j} \rightarrow\infty$ as $k \rightarrow
\infty$. Consider the second term
\[
\lim\sup_{k \rightarrow
\infty}\sum_{i=m}^{k-j-1}e^{-d(\overline{Y}{}^\alpha_{k-j}-\overline
{Y}{}^\alpha_i)}
\le\lim_{k \rightarrow\infty}
\sum_{l=1}^{k-j-m}e^{-dlw^\alpha(1-\varepsilon)}=\frac{1}{e^{d(1-\varepsilon
)w^\alpha}-1}.
\]
Similarly,
\[
\lim\inf_{k \rightarrow
\infty}\sum_{i=m}^{k-j-1}e^{-d(\overline{Y}{}^\alpha_{k-j}-\overline
{Y}{}^\alpha_i)}
\ge\frac{1}{e^{d(1+\varepsilon)w^\alpha}-1}.
\]
Hence,
%
\begin{equation}
\label{eqn74}
\lim_{k \rightarrow
\infty}T_{j,k}=\frac{1}{e^{dw^\alpha}-1} \qquad\mbox{for any fixed
$j$}.
\end{equation}
Substituting the results (\ref{eqn73}) and (\ref{eqn74}) into
(\ref{eqn72}) yields
%
\begin{eqnarray}
\label{eqn75} \lim_{k \rightarrow\infty}
\gamma(k)e^{d\overline{Y}{}^\alpha_{k-j}}&=&\frac{1}{({e^{dw^\alpha
j}-1})/({e^{dw^\alpha}-1})+{1}/({e^{dw^\alpha}-1})}\nonumber\\[-8pt]\\[-8pt]
&=&\frac{e^{dw^\alpha}-1}{e^{dw^{\alpha}j}}=\mu_j \qquad\mbox{for
fixed $j$}.\nonumber
\end{eqnarray}
Returning to (\ref{eqn71}), fix $n$.
%
\begin{eqnarray}
\label{eqn76}
&&-\sum_{j=1}^{k-1}\log\bigl[1+t\bigl(1+o_k(1)\bigr)\gamma(k)e^{d\overline
{Y}{}^\alpha_{k-j}} \bigr]\nonumber\\
&&\qquad=-\sum_{j=1}^{n-1}\log\bigl[1+t\bigl(1+o_k(1)\bigr)\gamma(k)e^{d\overline
{Y}{}^\alpha_{k-j}} \bigr]\\
&&\qquad\quad{}-\sum_{j=n}^{k-1}\log\bigl[1+t\bigl(1+o_k(1)\bigr)\gamma(k)e^{d\overline
{Y}{}^\alpha_{k-j}} \bigr].\nonumber
\end{eqnarray}
Equation (\ref{eqn75}) implies that each term in the sum of the
first expression on the right-hand side converges to
$\log(1+t\mu_j)$ as $k \rightarrow\infty$.
We need to show that $\lim_{k \rightarrow
\infty}\sum_{j=n}^{k-1}\log[1+t(1+o_k(1))\gamma(k)e^{d\overline
{Y}{}^\alpha_{k-j}} ]$
can be made arbitrarily small by choosing $n$ to be sufficiently
large (all terms in the sum are positive). Note that
$\gamma(k)e^{d\overline{Y}{}^\alpha_{k-j}} < \frac{1}{S_{j,k}}$. For
any $\varepsilon>0$ choose $n$ large enough so that
$\frac{\overline{Y}{}^\alpha_i}{i+1} \ge(1-\varepsilon)w^\alpha$ for all
$i \ge n$. For $j \ge n$,
\begin{eqnarray*}
S_{j,k}&=&\sum_{i=k-j}^{k-1}e^{d(\overline{Y}{}^\alpha_i-\overline
{Y}{}^\alpha_{k-j})}=\sum_{l=0}^{j-1}e^{d(\overline{Y}{}^\alpha
_{k-j+l}-\overline{Y}{}^\alpha_{k-j})}\\
& \ge&
\sum_{l=0}^{j-1}e^{dlw^\alpha(1-\varepsilon)}=\frac{e^{dw^\alpha
j(1-\varepsilon)}-1}{e^{dw^\alpha(1-\varepsilon)}-1}.
\end{eqnarray*}
Hence,
\[
\gamma(k)e^{d\overline{Y}_{k-j}} <
\frac{e^{dw^\alpha(1-\varepsilon)}-1}{e^{dw^\alpha
j(1-\varepsilon)}-1}<e^{-dw^\alpha(j-1)(1-\varepsilon)}.
\]
Choose $k$ large enough so that $o_k(1) < \varepsilon$. Then
\begin{eqnarray*}
&&
\lim_{k \rightarrow
\infty}\sum_{j=n}^{k-1}\log\bigl[1+t\bigl(1+o_k(1)\bigr)\gamma(k)e^{d\overline
{Y}{}^\alpha_{k-j}} \bigr]
\\
&&\qquad\le\lim_{k \rightarrow
\infty}\sum_{j=n}^{k-1}t(1+\varepsilon)e^{-dw^\alpha(j-1)(1-\varepsilon)}
\\
&&\qquad<\frac{t(1+\varepsilon)e^{-dw^\alpha(n-2)(1-\varepsilon)}}{e^{dw^\alpha
(1-\varepsilon)}-1}.
\end{eqnarray*}
Since the right-hand side goes to zero as $n \rightarrow\infty$ the
second term in the sum in (\ref{eqn76}) can be made arbitrarily
small by choosing $n$ sufficiently large.
\end{pf*}
\end{appendix}


%
\printaddresses

\mbox{}

\end{document}